\documentclass[review,3p]{elsarticle}

\usepackage{lineno}
\modulolinenumbers[5]
\usepackage{dsfont}
\usepackage{amsmath}
\usepackage{amsfonts}
\usepackage{amssymb}
\usepackage{multirow}
\usepackage{booktabs}
\usepackage[table,xcdraw]{xcolor}
\usepackage[T1]{fontenc}
\usepackage{epsfig}
\usepackage{amssymb}
\usepackage{color}
\usepackage[linesnumbered,ruled,vlined]{algorithm2e}
\usepackage{pdfpages}                                                 
\usepackage{graphicx}
\usepackage{tikz}
\usetikzlibrary{fit}
\usepackage{wrapfig}
\usepackage{epstopdf}
\usepackage[font=small]{caption} 
\usepackage[font=small,position=bottom]{subfig}
\usepackage{stfloats}
\usepackage{tabularx}
\usepackage{placeins}
\usepackage{mathrsfs}
\usepackage{enumerate}
\usepackage[normalem]{ulem}

\journal{Journal of \LaTeX\ Templates}

\SetKwProg{Fn}{Function}{}{End}
\SetKwFunction{Frec}{\textnormal{RecMRCM}}%
\SetAlgoLongEnd
\SetAlFnt{\textnormal}
\def\Om{\Omega}

\def\p{\partial}

\def\bu{{\bf u}}
\def\bx{{\bf x}}
\def\bn{{\bf n}}
\definecolor{btru}{rgb}{0.18,0.33,0.59}

\newcommand{\sizef}{\footnotesize}

\newcommand{\sizefonte}{\footnotesize}

\newcolumntype{Y}{>{\centering\arraybackslash}X}
\begin{document}

\begin{frontmatter}

%\title{Recursive formulation and parallel implementation of multiscale mixed methods for subsurface flows}

\title{Recursive formulation and parallel implementation of multiscale mixed methods}

%\tnotetext[mytitlenote]{Fully documented templates are available in the elsarticle package on \href{http://www.ctan.org/tex-archive/macros/latex/contrib/elsarticle}{CTAN}.}

%% Group authors per affiliation:

%\author{E. Abreu} \author{A. M. Esp\'irito Santo} \author{P. Ferraz} \author{F. Pereira} \author{L. G. C. Santos} \author{F. S. Sousa}
%
%\address{Instituto de Matem\'{a}tica, Estat\'{i}stica e Computa\c{c}\~{a}o Cient\'{i}fica, Universidade de Campinas, R. S\'{e}rgio Buarque de Holanda, 651, 13083-859, Campinas, SP , Brazil} 
%\address{Centro de Estudos do Petr\'{o}leo, Universidade de Campinas, R. Cora Coralina, 350, 13083-896 Campinas, SP, Brazil} 
%\address{Department of Mathematical Sciences, The University of Texas at Dallas, 800 W. Campbell Road, Richardson, TX 75080-3021, USA} 
%\address{Instituto de Ci\^{e}ncias Matem\'{a}ticas e de Computa\c{c}\~{a}o, Universidade de São Paulo, Av. Trabalhador S\~{a}o-carlense, 400, 13566-590, S\~{a}o Carlos, SP, Brazil}
%
%
\address[address1]{Instituto de Matem\'{a}tica, Estat\'{i}stica e Computa\c{c}\~{a}o Cient\'{i}fica, Universidade de Campinas, R. S\'{e}rgio Buarque de Holanda, 651, 13083-859, Campinas, SP , Brazil}
\address[address2]{Centro de Estudos do Petr\'{o}leo, Universidade de Campinas, R. Cora Coralina, 350, 13083-896 Campinas, SP, Brazil}
\address[address3]{Departamento de Matemática Pura e Aplicada, 
Universidade Federal do Rio Grande do Sul, Av. Bento Gonçalves, 9500, 13083859, Porto Alegre, RS, Brazil}
\address[address4]{Department of Mathematical Sciences, The University of Texas at Dallas, 800 W. Campbell Road, Richardson, TX 75080-3021, USA}
\address[address5]{Instituto de Ci\^{e}ncias Matem\'{a}ticas e de
Computa\c{c}\~{a}o, Universidade de S\~ao Paulo, Av. Trabalhador S\~{a}o-carlense, 400, 13566-590, S\~{a}o Carlos, SP, Brazil}

%% or include affiliations in footnotes:
%\author[mymainaddress,mysecondaryaddress]{Elsevier Inc}
%\ead[url]{email@email.com.br}

\author[address1]{E. Abreu}
\author[address2]{P. Ferraz\corref{mycorrespondingauthor}}
\author[address3]{A. M. Esp\'irito Santo}
\author[address4]{F. Pereira}
\author[address2]{L. G. C. Santos}
\author[address5]{F. S. Sousa}
\cortext[mycorrespondingauthor]{Corresponding author \\ \textit{ \indent E-mail
Address:} paola.ferraz@gmail.com, pferraz@ime.unicamp.br}
%\ead{paola.ferraz@gmail.com}

\begin{abstract}

Multiscale methods for second order elliptic equations based on
non-overlapping domain decomposition schemes have great potential to take
advantage of multi-core, state-of-the-art parallel computers.  These methods
typically involve solving local boundary value problems followed by the
solution of a global interface problem. Known iterative procedures 
for the solution of the interface problem have typically slow
convergence, increasing the overall cost of the multiscale solver.  To overcome
this problem we develop a scalable recursive solution method
for such interface problem that replaces the global problem by a family of
small interface systems associated with adjacent subdomains, in a hierarchy of
nested subdomains. Then, we propose a novel parallel algorithm to implement
our recursive formulation in multi-core devices using the Multiscale Robin
Coupled Method by Guiraldello et al. (2018), that can be seen as a
generalization of several multiscale mixed methods.  Through several numerical
studies we show that the new algorithm is very fast and exhibits excellent strong and weak 
scalability.  We consider very large problems, that can have billions of
discretization cells, motivated by the numerical simulation of subsurface
flows.
\end{abstract}


\end{frontmatter}

%\linenumbers

% Introduction:
%\input{section/section_1}
\section{Introduction}\label{sec:3.1}

Multiscale methods have been developed in the last few decades to approximate
efficiently problems involving second order elliptic partial differential
equations.  These problems are very important in several areas of research, in
particular in applications to oil reservoir simulation with high contrast in
heterogeneity. Despite considerable advances in computational processing
capability and storage, traditional methods that have been used in mainstream
oil reservoir simulators are not capable of dealing with problems involving
billions of elements in the discretization of large computational regions.
Very large reservoirs of interest to the industry can be found, for instance,
in the Brazilian pre-salt layer and acceptable accuracy in numerical
simulations require a considerable number of elements. 
As more variables and processes are taken into account 
to accurately resolve fine scale details of real life models,  
resource efficiency is an important requirement. 
A number of multiscale
methods have been developed to overcome the computational challenges posed by
these simulations and ensure acceptable precision of numerical solutions.
Domain decomposition techniques divide the global
domain into subregions that may be overlapping or non-overlapping,
facilitating the use of parallelization techniques. Local
solutions, called multiscale basis functions, are constructed through solutions of
boundary value problems within each subdomain. These functions
retain fine mesh information and are employed as building blocks to construct
global approximations for the problem at hand. The key idea is to obtain an
approximate solution considering unknowns defined on a coarse scale, 
and thus reducing drastically
the number of unknowns with respect to the fine mesh. The multiscale basis
functions are then used to reconstruct the fine scale solution from the coarse problem. 

Two major classes of multiscale methods can be identified: methods in the
context of finite elements such as the Multiscale Finite Element 
Methods (MSFE) \cite{hou97,efmults} and the Generalized Multiscale 
Finite Element Method (GMsFEM) \cite{EGH13},
and those that use finite volume such as the Multiscale Finite Volume Methods 
(MSFV) 
\cite{jenny2003,jenny2005,Lee2008,Lunati2006,Jenny2009,Sokolova2019}. 
On the other hand, extensions of these multiscale methods were formulated
to be used as preconditioners in iterative algebraic solvers     
\cite{H08,HJ11,manea16,manea17,manea19}.
The formulation of multiscale methods are frequently naturally parallelizable
and some methods were implemented in multi-core CPU/GPU systems 
(see, \cite{manea16,manea17,manea19,ZYSW17,PED18}). 
The largest three-dimensional problem considered in these references 
has $128$ million discretization
cells, and was run in CPU/GPU clusters.
In our work the focus is on the family of multiscale mixed methods 
composed by the Multiscale Mortar Mixed Finite Element Method
(MMMFEM) \cite{APWY07,ganis2009,wheeler2012,Ahmed2019}, the Multiscale
Hybrid-Mixed Method (MHM) \cite{Valentin2013,Harder2013,Araya2017}, the
Multiscale Mixed Method (MuMM) \cite{Pereira2014,APGK19} and the Multiscale
Robin Coupled Method (MRCM) \cite{guiraldello2018multiscale,Guiraldello2019}
that has been more recently introduced in the literature.  
For these methods, the coarse scale is defined by the skeleton
of an underlying domain decomposition where the subdomains are coupled using
distinct interface conditions.  
The MMMFEM couples subdomains through a continuous pressure and 
weak continuity of normal fluxes.
Thus, a post-processing step is inevitable to produce velocity fields with
continuous normal components on the fine grid.  On the other hand the MHM
couples subdomains through the imposition of continuous normal flux components,
and the pressure is weakly continuous. The MuMM is a multiscale domain
decomposition method based on the work of \cite{DPLRW93} where the Robin boundary
conditions are used to obtain local solutions. In the MuMM the continuity of
normal component of fluxes as well as the pressure are weakly imposed.
Finally, there is the MRCM that also utilizes the Robin coupling conditions between
subdomains and generalizes the above mentioned multiscale mixed approaches.  In
\cite{guiraldello2018multiscale} it is shown that the MMMFEM and MHM can be
seen as members of a family of multiscale methods parametrized by the Robin
condition coefficient. The MuMM can also be seen as a particular case of the
MRCM, when considering piecewise constant spaces set at the skeleton of the
decomposition.

Our contribution in this work is twofold. 
First we introduce a recursive formulation for a family of
multiscale mixed methods that is used to construct a new interface solver developed
specifically for parallel processing in multi-core systems.
The new recursive formulation can be seen as a
variational formulation of the procedure recently
introduced (and referred to as a multiscale direct solver) 
in \cite{APGK19}.  Then, we propose a novel parallel algorithm  based on
the recursive formulation. Through a careful analysis for large problems 
we show that the proposed algorithm is very fast
and exhibits excellent scalability, both strong and weak. 
We consider larger problems as well as larger number of processing cores than in
existing parallel results produced by multiscale methods for elliptic equations.
For more details of the new recursive approach see \cite{Paolaphd2019}.

This work is organized as follows. In Section \ref{sec:3.2} we briefly review
the MRCM method.  In Section \ref{sec:3.3} we describe in details the recursive
formulation and its parallel implementation. We discuss the connection between
the MuMM and the MRCM in Section \ref{sec:3.4} and in Section \ref{sec:3.5} we
present numerical experiments to show the excellent scalability of our
proposed method. In Section 6 we discuss our work with other
parallel implementations. 
Finally, in Section 7 we present our concluding remarks.

% Review of MRCM:
%\input{section/section_2}
\section{A review of the Multiscale Robin Coupled Method}\label{sec:3.2}

The Multiscale Robin Coupled Method (MRCM) introduced in
\cite{guiraldello2018multiscale} is a multiscale mixed method based on a
non-overlapping domain decomposition where subdomains are coupled through
weak continuity of pressure and normal across the interfaces
between subdomains. The parameter appearing in the Robin condition used in the
local boundary value problems associated with the subdomains
determines the relative importance of Dirichlet or Neumann boundary condition
in the coupling of subdomains. The result is that for small (resp. large)
values of this parameter, the solution produced by the MRCM converges to the
solution of the MMMFEM (resp. MHM), a property that is well illustrated and
explored in \cite{guiraldello2018multiscale}. This parameter plays an important
role in the approximation of two-phase flows in high-contrast porous media, as
can be seen in \cite{Fran2019}. Another aspect of this method is that
it introduces great flexibility in the choice of interface spaces 
for normal fluxes and pressures at the skeleton of the decomposition
(see \cite{Guiraldello2019}). It is also observed in
\cite{guiraldello2018multiscale} that the variational formulation of the MRCM
is an extension of the MuMM, that was originally introduced
as an iterative method, and can be recovered by
a suitable choice of parameters for the MRCM. 

In this section, we recall the key aspects of the MRCM. To briefly introduce
the variational formulation, consider a rectangular domain $\Omega \subset
\mathbb{R}^{d}, d \in \{2,3\}$, with a  Lipschitz boundary $\partial\Omega$,
defined for the following pressure-velocity problem in mixed form,
\begin{align}
\nabla \cdot {\bf u} &= f({\bf x}),
\qquad
{\bf u} = -K({\bf x}) \nabla p({\bf x}), 
\qquad \bx \in \Omega,  \label{eq:ellip.1} \\
p &= g_D, \qquad  \bx \in \partial\Omega_D, \label{eq:ellip.2} \\
{\bf u} \cdot \check{\bf n} &= g_N, \qquad \bx \in \partial\Omega_N, \label{eq:ellip.3}
\end{align}
where $\bu = \bu(\bx)$ is the Darcy's flux and $p(\bx)$ is the fluid pressure. 
The absolute permeability is given by $K(\bx)$, a symmetric positive definite 
tensor, and $\check{\bn}$ is the unit outward normal vector to $\p\Omega$. 

The domain decomposition formulation of the MRCM is performed directly in the
discrete form of the system \eqref{eq:ellip.1}-\eqref{eq:ellip.3}. Thus we
start by decomposing the domain $\Om$ into $m$ non-overlapping subdomains
$\Omega^i \, ,\, i = 1,\dots,m$, with reference size $H$, where
\begin{equation}\label{eq:decomp.1}
\Omega = \bigcup\limits_{i=1}^{m} {\Omega}^i,  
\quad
{\Omega}^k \cap {\Omega}^i = \emptyset, 
\quad
i \neq k,
\end{equation}
each with a well-defined Lipschitz boundary $\partial\Omega^i$. Let $\Gamma =
\cup_i \partial\Omega^i \setminus \partial\Omega$, be the skeleton of the
domain decomposition, and 
\begin{equation}\label{eq:decomp.2}
\Gamma^i = \Gamma \cap \partial\Omega^i,
\quad
\Gamma^{ik} = \Gamma^{ki} = 
\partial\Omega^i \cap \partial\Omega^k.
\end{equation}
We refer to $\Gamma^{ik} = \Gamma^{ki}$ as the interface between the subdomains 
$\Omega^i$ and $\Omega^k$.
Additionally, let us define two types of normal vectors. One denoted by
$\check{\bn}^i$ is simply the normal vector pointing outward of subdomain
$\Om^i$. The second, denoted as $\check{\bn}$ with no superscript, will have a
global definition on $\Gamma$, that is for every $\Gamma^{ik}\subset \Gamma$,
it points towards the subdomain with maximum index value ($\max\{i,k\}$). This
will be used as a reference vector in the variational formulation, to uniquely
identify the direction of fluxes over each interface of $\Gamma$.

\begin{figure}[h!]
\centering
\includegraphics[width=1.0\textwidth]{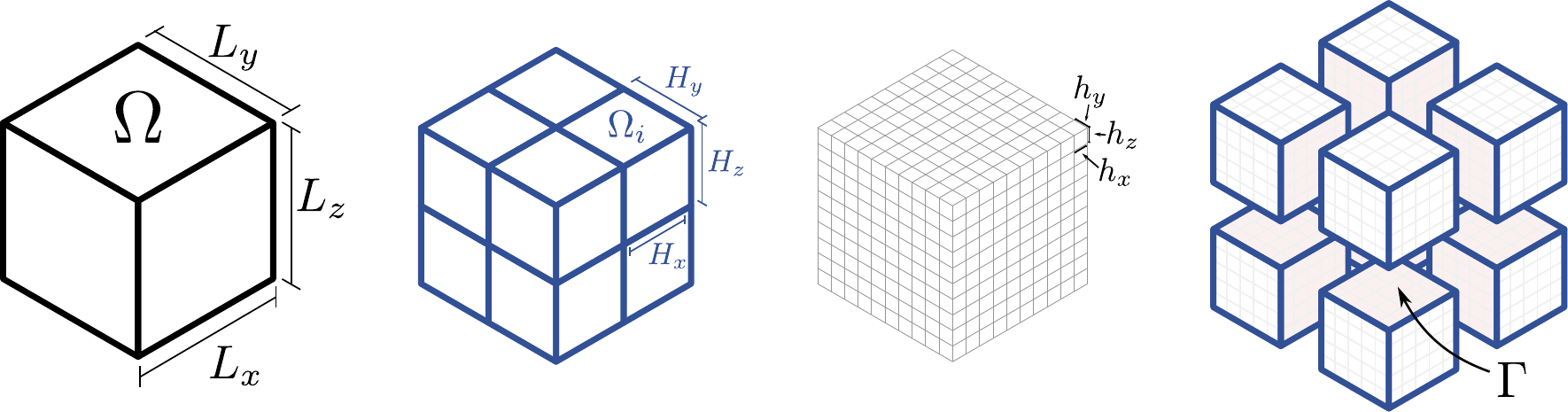}
\caption{
Representations of a three-dimensional domain decomposition of $\Om$. 
On the leftmost image, the complete domain is shown with sizes $L_x, L_y$ and
$L_z$. The second and the third pictures show the coarse scale ($H$) and the
fine scale ($h$), respectively. Rightmost picture depicts $\Gamma$, the
skeleton of the decomposition and is composed by subdomain interfaces.  }
\label{fig:hbarscale}
\end{figure}

Let $\mathcal{T}_h^i$ be a regular mesh discretizing $\Om^i$, with reference
size $h \ll H$ (see Figure \ref{fig:hbarscale}), where it is possible to define
the lowest order Raviart-Thomas spaces for velocity and pressure, say ${\bf
V}_{h}^i \subset H(\mbox{div},\Omega^i)$ and $Q_h^i \subset L^2(\Omega^i)$,
respectively (their definition can be seen in
\cite{guiraldello2018multiscale}). We will also need the vector space ${\bf
V}_{h,g_N}^i \subset {\bf V}_{h}^i$ of the functions in ${\bf V}_{h}^i$
satisfying the Neumann boundary conditions in \eqref{eq:ellip.3}. 

Finally, the variational formulation of the MRCM introduces unknowns $U_H$ and
$P_H$ that are defined only on the skeleton $\Gamma$ of the domain
decomposition. For that purpose, interface spaces are needed, that are 
defined as subspaces of the set of piecewise constant functions 
\begin{equation}
F_h(\mathcal{E}_h) = \{ f:{\mathcal{E}}_h\to
\mathbb{R}~|~f|_e\,\in\,\mathbb{P}_0~, ~\forall\,e\,\in\,{\mathcal{E}}_h \},
\label{eq:spaces.edges}  
\end{equation} 
where $\mathcal{E}_h$ is the set of
all edges/faces of $\Gamma$. Hence, we can take $\mathcal {P}_h =
F_h(\mathcal{E}_h)$ as the pressure interface spaces, and $\mathcal{U}_h$ as
being the subspace of $F_h(\mathcal{E}_h)$ of the functions that are zero when
$\beta^i$, 
the Robin condition parameter defined as a function on $F_h(\mathcal{E}_h)$,
vanishes on both sides of the interface. The multiscale formulation of
the MRCM is defined over the coarse subspaces $\mathcal {P}_H \subset \mathcal
{P}_h$ and $\mathcal {U}_H \subset \mathcal {U}_h$, and formalized below:

\paragraph{Discrete variational formulation of the MRCM} 
Find the local solution $({\bf u}^i_h,p^i_h)\,\in\,{\bf V}_{h,g_{N}}^i\times Q_h^i$, for
$i=1,\ldots,m$, and
$(U_H,P_H)\,\in\, \mathcal{U}_H \times \mathcal{P}_H$
such that
\begin{align}
(K^{-1}{\bf u}^i_h,{\bf v})_{\Omega^i}-(p^i_h,\nabla\cdot {\bf v})_{\Omega^i}
+(P_H - \beta^i U_H \check{\bf n}^i \cdot \check{\bf n}+\beta^i\,{\bf u}^i_h \cdot \check{\bf
n}^i,{\bf v} \cdot \check{\bf n}^i)_{\Gamma^i} &=
-(g_D\,,{\bf v} \cdot \check{\bf n}^i)_{\partial\Omega^i\cap\partial\Omega_D}~,
\label{eq11.2} \\
(q,\nabla\cdot {\bf u}^i_h)_{\Omega^i} &= (f,q)_{\Omega^i}~, \label{eq12.2}
\end{align}
with the following interface conditions
\begin{align}
\sum_{i=1}^m ({\bf u}^i_h \cdot \check{\bf n}^i,\,M_H)_{\Gamma^i} &=  0~,
\label{eq13.2}\\
\sum_{i=1}^m (\beta^i\,({\bf u}^i_h \cdot \check{\bf n}^i - U_H \check{\bf n} \cdot \check{\bf n}^i),\,
V_H \check{\bf n}^i \cdot \check{\bf n})_{\Gamma^i} &= 0~, \label{eq14.2}
\end{align}
hold for all $({\bf v},q)\,\in\,
{\bf V}_{h,0}^i\times Q_h^i~,~\forall\,i=1,\ldots,m$, and for all
$(V_H,M_H)\,\in\, \mathcal{U}_H\times \mathcal{P}_H$.\\

More details about this variational formulation, as well as the well-posedness
of the discrete system, can be seen in \cite{guiraldello2018multiscale}.
The final global solution $(\bu_h,p_h)$ of
\eqref{eq:ellip.1}-\eqref{eq:ellip.3} is written as a combination of
the local solutions  $(\bu^i_h,p^i_h)$.

\subsection{Mixed multiscale basis functions}\label{sec:mmbf}

An efficient implementation of mixed multiscale methods can be achieved by
writing the final solution in terms of a set of mixed multiscale basis
functions (hereafter referred as MMBF's), a procedure already discussed by
other authors, such as Ganis \& Yotov \cite{ganis2009}, Francisco et al.
\cite{Pereira2014} and more recently by Guiraldello et al.
\cite{guiraldello2018multiscale}. Following the ideas already presented by
these authors, especially the later, we recall this procedure to introduce
the notation for our recursive formulation of the MRCM. 

We start with an additive  decomposition of the local solutions $({\bf
u}^i_h,p^i_h)$ in $\Omega^i$, as 
\begin{align}
{\bf u}^i_h & =  \widehat{\bf u}^i_h + \bar{\bf u}^i_h, \label{eq:addu}\\
p^i_h & =  \widehat{p}^i_h + \bar{p}^i_h, \label{eq:addp}
\end{align}
where $(\widehat{\bf u}^i_h,\widehat{p}^i_h) \in {\bf V}_{h,0}^i\times Q_h^i$
represents the homogeneous part, i.e., the solution of the local problem
\eqref{eq11.2}-\eqref{eq12.2} with given Robin boundary conditions (given $U_H$
and $P_H$), and  vanishing source and external boundary data, while  $(\bar{\bf
u}^i_h,\bar{p}^i_h) \in {\bf V}_{h,g_N}^i\times Q_h^i$ is the solution of the
local problem \eqref{eq11.2}-\eqref{eq12.2} with vanishing Robin boundary
conditions ($U_H=P_H=0$), nonzero source, and external boundary data.

The solution of the homogeneous part $(\widehat{\bf u}^i_h,\widehat{p}^i_h)$
can be obtained as a linear combination of MMBF's, which can be constructed by
properly setting $U_H$ and $P_H$. Consider $\{\phi^{j}\}_{1 \leq j \leq
n_U}$ and $\{\psi^{j}\}_{1 \leq j \leq n_P}$ a finite element basis for
the coarse interface spaces $\mathcal{U}_H$ and $\mathcal{P}_H$, respectively,
where $n_U = \dim(\mathcal{U}_H)$ and $n_P = \dim(\mathcal{P}_H)$.  Then,
the interface variables $U_H$ and $P_H$ can be written as
\begin{equation}\label{eq:coarse.space.var}
U_H = \sum\limits_{{j}=1}^{n_U} X_{j} \phi^{j},
\quad
P_H = \sum\limits_{{j}=1}^{n_P} X_{{j}+n_U} \psi^{j},
\end{equation}
where the coefficients $X = (X_1, \dots, X_{n})^T$ are to be determined
later. Define $\mathcal{J}$ as the set of global indices 
of the interface degrees
of freedom, such that $|\mathcal{J}| = n = n_U + n_P$.
Also define $\mathcal{J}^{i}$ as the set of interface degrees of freedom 
associated with $\Omega^i$
whose support is on the boundary $\Gamma^i$, such that 
$|\mathcal{J}^{i}| = n^i$.    
For every $j \in \mathcal{J}^{i}$, the multiscale basis function in $\Omega^i$,
denoted here as $\{\boldsymbol{\Phi}^i_{k_j}, \Psi^i_{k_j} \}_{1\le k_j \le
n^i}$, 
are given by the following local problems:
\begin{itemize}
\item
If $1 \le j \le n_U $, solve problem \eqref{eq11.2}-\eqref{eq12.2} with
boundary data $U_H = \phi^j$, $P_H = 0$:

Find $(\boldsymbol{\Phi}^i_{k_j}, \Psi^i_{k_j})\,\in\,{\bf V}_{h,0}^i\times Q_h^i$, 
such that
\begin{align}
(K^{-1}\boldsymbol{\Phi}^i_{k_j},{\bf v})_{\Omega^i}-
(\Psi^i_{k_j},\nabla\cdot {\bf v})_{\Omega^i} +
(\beta^i\,\boldsymbol{\Phi}^i_{k_j}\cdot\check{\bf n}^i,{\bf v}\cdot
\check{\bf n}^i)_{\Gamma^i}
&= (\beta^i \phi^j \check{\bf n}^i\cdot\check{\bf n},{\bf v}\cdot\check{\bf n}^i)_{\Gamma^i}, 
\label{eq11Hhatbasisa} \\ 
(q,\nabla\cdot \boldsymbol{\Phi}^i_{k_j})_{\Omega^i} &= 0,
\label{eq12Hhatbasisa}
\end{align}
hold for all $({\bf v},q)\,\in\,{\bf V}_{h,0}^i\times Q_h^i$. 

\item
If $n_U < j \le n$,  solve problem \eqref{eq11.2}-\eqref{eq12.2} with boundary data $U_H = 0$,
$P_H = \psi^{j-n_U}$: 

Find $(\boldsymbol{\Phi}^i_{k_j}, \Psi^i_{k_j})\,\in\,{\bf V}_{h,0}^i\times Q_h^i$, 
such that
\begin{align}
(K^{-1}\boldsymbol{\Phi}^i_{k_j},{\bf v})_{\Omega^i}-(\Psi^i_{k_j},\nabla\cdot {\bf v})_{\Omega^i} +
(\beta^i\,\boldsymbol{\Phi}^i_{k_j}\cdot\check{\bf n}^i,{\bf v}\cdot\check{\bf
n}^i)_{\Gamma^i}
&=
-(\psi^{j-n_U},{\bf v}\cdot\check{\bf n}^i)_{\Gamma^i}, \label{eq11Hhatbasisb} \\ 
(q,\nabla\cdot \boldsymbol{\Phi}^i_{k_j})_{\Omega^i}&= 0, \label{eq12Hhatbasisb}
\end{align}
hold for all $({\bf v},q)\,\in\,{\bf V}_{h,0}^i\times Q_h^i$.
\end{itemize}

In the variational formulations above, the functions $\phi^j$ and $\psi^{j}$
depend on the interface space considered. An exploration of several choices for
interface spaces, both polynomial and informed spaces, are considered in
\cite{Guiraldello2019}.  The homogeneous local solutions $(\widehat{\bf
u}^i_h,\widehat{p}^i_h)$ are then written as a linear combination of the
multiscale basis functions, $\{\boldsymbol{\Phi}^i_{k_j}, \Psi^i_{k_j} \}_{1\le
k_j \le n^i}$, as
\begin{equation}
\widehat{\bf u}^i_h = \sum_{j \in \mathcal{J}^{i}} {X_{j}
\boldsymbol{\Phi}^i_{k_j}},
\qquad
\widehat{p}^i_h = \sum_{j \in \mathcal{J}^{i}} {X_{j} \Psi^i_{k_j}}.
\label{eq:uhatphat}
\end{equation}

The local problems \eqref{eq11Hhatbasisa}-\eqref{eq12Hhatbasisb} can be solved
by any discretization that delivers both pressure and normal fluxes at the
skeleton $\Gamma$ of the decomposition. In \cite{guiraldello2018multiscale},
the authors perform a conservative finite volume discretization, while in this
work, we use the (equivalent) lowest order Raviart-Thomas (RT0) spaces
for the interface unknowns, such as in
\cite{DPLRW93,raviart1977mixed,EA14,RT91,Paolaphd2019}. Although conveniently
parallelizable, given the local nature of the problems involved, the
computation of a large set of MMBF's can still be very expensive, even in
multi-core high-performance computers. 

\subsection{Interface system}

The use of multiscale basis functions allows us to build a linear system for the
interface unknowns alone \cite{guiraldello2018multiscale,ganis2009,APGK19}.
The procedure consists of substituting the solution \eqref{eq:addu}-\eqref{eq:addp}
written as a linear combination of the MMBF's (as in \eqref{eq:uhatphat}) in
the coarse scale continuity conditions \eqref{eq13.2}-\eqref{eq14.2}. The next
step is to substitute the interface unknowns by the linear combinations in
\eqref{eq:coarse.space.var} and test $V_H$ and $M_H$ appearing in
\eqref{eq13.2}-\eqref{eq14.2} for all basis functions spanning $\mathcal{U}_H$
and $\mathcal{P}_H$. As a result, we end up with a linear system of the form 
\begin{equation}
\mathbf{A}\, X = \mathbf{b}, \label{eq:linsys}
\end{equation}
where the unknown vector $X = (X_1,\ldots,X_{n})^{T}$ is formed by the
coefficients of the linear combinations in \eqref{eq:coarse.space.var}. The
entries of matrix $\mathbf{A}$ are, for $j = 1,\dots, n$
\begin{equation}
 a_{rj} = \left.
  \begin{cases}
  \sum_{i=1}^m \left(\beta^i\,(\boldsymbol{\Phi}^i_{k_j}\cdot \check{\bf n}^i-\varphi^{j} \check{\bf n}^i\cdot\check{\bf n}),
  \phi^{{r}}\, \check{\bf n}^i\cdot\check{\bf n}\right)_{\Gamma^i}, &
  \text{for } 1 \leq {r} \leq n_U \\
  \sum_{i=1}^m{\left( \boldsymbol{\Phi}^i_{k_j} \cdot \check{\bf n}^i , \psi^{{r}}
  \right)_{\Gamma^i}}, & 
  \text{for } n_U < {r} \leq n
  \end{cases}
  \right.\label{eq:brhs}
\end{equation}
where $\varphi^j = \phi^j$ if $1 \le j \le n_U$ and zero otherwise. 
As for the right hand side vector $\mathbf{b}$, 
computing its entries involves the particular solutions $\bar{\bu}^i_h$, yielding
\begin{equation}
 \mathbf{b}_{r} = \left.
  \begin{cases}
    -\sum_{i=1}^m (\beta^i\,(\bar{\bu}^i_h \cdot \check{\bf n}^i), \phi^{{r}} \,
    \check{\bf n}^i\cdot \check{\bf n})_{\Gamma^i}, 
    & \text{for } 1 \leq {r} \leq n_U \\
    -\sum_{i=1}^m {\left( \bar{\bu}^i_h \cdot \check{\bf n}^i , \psi^{{r}}
    \right)_{\Gamma^i}}
    , & \text{for } n_U < {r} \leq n. 
  \end{cases}
  \right.
\end{equation}

Lastly, the local final solution $({\bf u}^i_h,p^i_h)$ in $\Omega^i$, given by
\eqref{eq:addu}-\eqref{eq:addp}, can be written as
\begin{equation}
{\bf u}^i_h = \sum_{j \in \mathcal{J}^{i}} {X_{j} \boldsymbol{\Phi}^i_{k_j}}
+ \bar{\bf u}^i_h,
\qquad
{p}^i_h = \sum_{j \in \mathcal{J}^{i}} {X_{j} \Psi^i_{k_j}}
+ \bar{p}^i_h.
\label{eq:final}
\end{equation}

Although quite efficient due to the reduced number of unknowns, this procedure
still needs the global assembly and resolution of the  non-symmetric linear
system \eqref{eq:linsys}, that, if not properly done, can hinder the parallel
efficiency of the overall method. In the following sections, we will introduce
a new naturally parallelizable methodology to localize and decompose the
interface problems for maximum efficiency.

% Recursive interface approach:
%\input{section/section_3}
\section{Recursive formulation}\label{sec:3.3}

We define the recursive formulation for the MRCM in terms of a hierarchy of 
nested decompositions of the domain $\Omega$ where the MRCM is applied recursively.  
The proposed method approximates the solution of the global problem by 
the solution of a family of smaller problems that fit well into  multi-core 
parallel machines (see also \cite{Paolaphd2019}). 
The general idea is to start by using the MRCM on a two-subdomain decomposition  
on $\Omega$, where each subdomain is 
successively decomposed in two smaller adjacent subdomains until a 
last stage is reached. The global interface problem is then
replaced by a family of small interface systems.
For simplicity, in this discussion we assume $\Omega$ to be a parallelepiped 
and all subdomains are cubes.

\subsection{A hierarchy of decompositions of the domain $\Omega$}\label{sec:notation}

Let us introduce the notation. 
We define a hierarchy of domain
decompositions in level $\ell$ given by,
\begin{equation}\label{eq:dd.l}
\Omega = \bigcup\limits_{i=1}^{m^\ell}~\Omega^{i,\ell},
\quad
m^\ell = 2^\ell,
\quad
\ell=0,\ldots,\mathcal{L}.
\end{equation}
such that in the zero-th level there is no decomposition, i.e., 
$\Omega^{1,0} = \Omega$.
The subdomains of the finest decomposition have sides of size $H$.
To define the hierarchy of decompositions of $\Omega$ we define each subdomain 
of level $\ell$ as being 
composed by the union of two subdomains of the decomposition of $\Omega$ on
level $\ell + 1$, 
\begin{equation}\label{eq:dd.2}
\Omega^{i,\ell} = \Omega^{2i-1,\ell+1}\,\cup\,\Omega^{2i,\ell+1},
\quad
i = 1,\ldots,m^{\ell}.
\end{equation}
For each level $\ell$ we define $\Gamma^{\cdot,\ell} = \cup_{i=1}^{m^{\ell}} 
\p\Omega^{i,\ell} \setminus \p\Omega$, 
as the skeleton of its associated domain decomposition where the ``$\cdot$''
superscript is to differentiate when the skeleton is 
defined on levels. 
We set subdomain interface for each level  as 
$\Gamma^{i,\ell} = \Gamma^{\cdot,\ell} \cap \p\Omega^{i,\ell}$ (for $\ell = 0$ we
have $\Gamma^{\cdot,0} = \emptyset$ by definition) and set $\Gamma^{ik,\ell} =
\Gamma^{ki,\ell} = \Omega^{i,\ell} \cap \Omega^{k,\ell}$ as the interface
between two subdomains on level $\ell$. 
Also set
\begin{equation}\label{eq:interface.1}
\gamma^{i,\ell} = \p\Omega^{2i-i,\ell+1}\, \cap \,\p\Omega^{2i,\ell+1},
\quad
i = 1,\ldots,m^{\ell},
\end{equation}
as the interface between two subdomains on level $\ell+1$ that compose $\Omega^{i,\ell}$
on level $\ell$, such that we are able to write the skeleton 
of the decomposition on each level as 
\begin{equation}\label{dd.gamma.i}
\Gamma^{\cdot,\ell+1} = \Gamma^{\cdot,\ell} \cup 
\left(\bigcup\limits_{i=1}^{m^\ell} \gamma^{i,\ell}\right),
\quad
\ell=0,\ldots,\mathcal{L}.
\end{equation}

\begin{figure}[h!]
	\centering
	\includegraphics[width=1\textwidth]{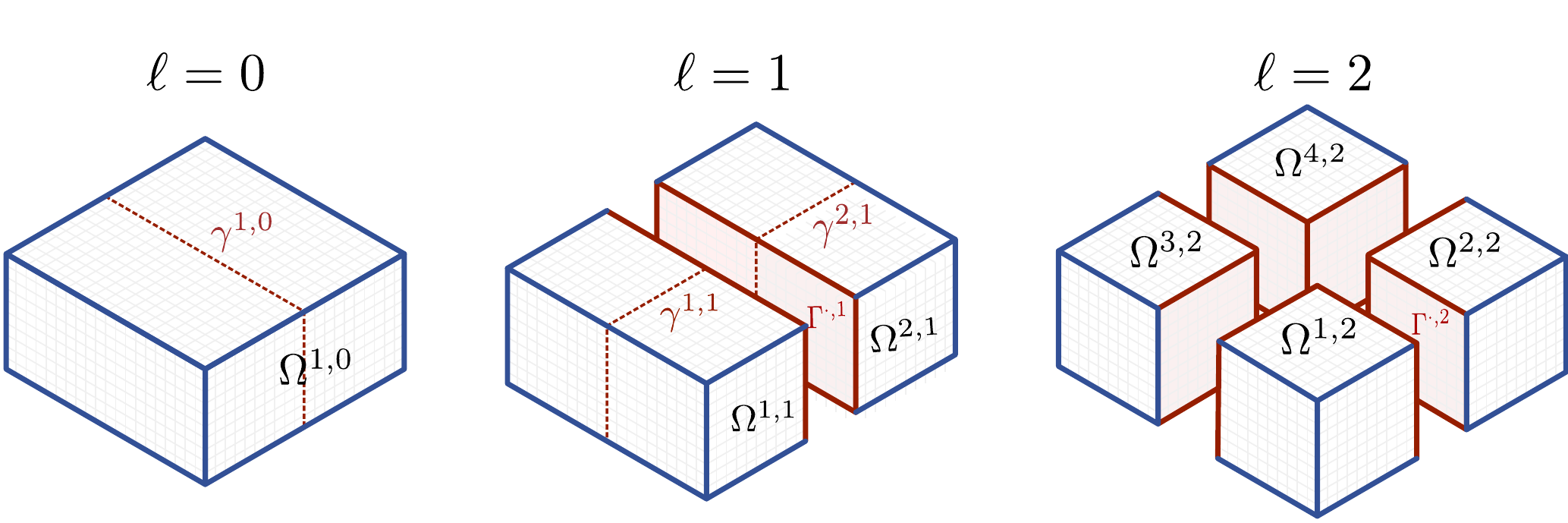}
	\caption{
		Representation of a sequence of domain decompositions form left to
        right.
        We begin at level 0 with the whole domain, where we will
        perform a Recursive MRCM step over $\gamma^{1,0}$, decomposing
        the domain into two subdomains. At level 1, we have
        the subdomains $\Omega^{1,1}$ and $\Omega^{2,1}$, with skeleton
        $\Gamma^{\cdot,1}$.
        At this level, we perform two steps of the Recursive MRCM on each 
        subdomain over $\gamma^{1,1}$ and $\gamma^{2,1}$ with each subdomain
        decomposed into two new subdomains. 
        We reach the finest level 2, where we
        have the \textit{finest} subdomain mesh with four subdomains and 
        skeleton of the decomposition $\Gamma^{\cdot,2}$.  
        }
	\label{fig:dd.6}
\end{figure}

For the interface spaces $\mathcal{U}_H$ and $\mathcal{P}_H$, we consider a 
finite element basis functions $\{\phi^{j}\}_{1 \leq j \leq n_U}$ and
$\{\psi^{j}\}_{1 \leq j \leq n_P}$ on the skeleton of the
\textit{finest} decomposition $\Gamma^{\cdot,\mathcal{L}}$ such that 
they have support on faces with size $H \times H$.
%of the subdomains of the \textit{finest} partition (level $\mathcal{L}$). 
In the recursive formulation, we define $\mathcal{J}^{\cdot,\ell}$ as the total set of
indices of interface degrees of freedom on level $\ell$, such that 
\begin{equation}\label{eq:degrees_freedom}
\mathcal{J}^{\cdot,0} \subset \ldots \subset \mathcal{J}^{\cdot,\ell} \subset \ldots \subset 
\mathcal{J}^{\cdot,\mathcal{L}}.
\end{equation}
We also define $\mathcal{J}^{i,\ell}$ as the set of interface degrees of freedom 
associated with $\Omega^{i,\ell}$ whose support is on the boundary $\Gamma^{i,\ell}$, 
such that $|\mathcal{J}^{i,\ell}| = n^{i,\ell}$.
Lastly, define $\xi^{i,\ell}$ as the interface degrees of freedom whose
support is on $\gamma^{i,\ell}$. 
Figure \ref{fig:dd.6} shows a two-level domain
decomposition sequence and its interfaces. 
Now we are ready to define the recursive formulation of
\eqref{eq11.2}-\eqref{eq14.2} to 
find the approximate solution $(\bu_h,p_h)$ of
\eqref{eq:ellip.1}-\eqref{eq:ellip.3}.

%\pagebreak

\subsection{Recursive formulation}\label{sec:idea}

The recursive formulation consists of the following steps.
Approximate \eqref{eq:ellip.1}-\eqref{eq:ellip.3}
by the MRCM where the domain $\Omega$ is decomposed in two subdomains. 
Within this decomposition, a family of MMBFs has to be computed for each subdomain
$\Omega^{i,1}, i = 1,2$. For each subdomain of level $\ell$, $\ell \geq 2$, we
follow the same procedure within $\Omega^{i,\ell}$, $i=1,\ldots,m^{\ell}$
subdomains.
Then, the MMBFs are computed by the use of the MRCM restricted to
each $\Omega^{i,\ell}$, $i=1,\ldots,m^{\ell}$.
This is achieved by decomposing $\Omega^{i,\ell}$ into two smaller subdomains and
following the usual steps of the MRCM for a two-subdomain decomposition.
We proceed from coarser ($\ell = 0$) to finer decompositions 
($\ell = \mathcal{L}$) by approximating the local problems by the MRCM, until 
the \textit{finest} decomposition is reached. 
At this point in the formulation we approximate the
solution of the MMBFs using a mixed finite element method (MFEM) or equivalently,
a finite volume method (FVM).
We define the recursive formulation of the MRCM in terms of a hierarchy of 
nested decompositions of $\Omega$ where the MRCM is applied recursively.  
We refer to this formulation as the Recursive MRCM, and it is introduced  
in Algorithm \ref{alg:recursive}.
Next we discuss Algorithm \ref{alg:recursive} in detail. 

\begin{algorithm}[ht!]
\DontPrintSemicolon

\Fn{\Frec{$\Omega^{i,\ell}$}} 
{
    \If{$\ell = \mathcal{L}$}
    {
   	
        Compute
        $\{\boldsymbol{\Phi}^{i,\mathcal{L}}_{s},\Psi^{i,\mathcal{L}}_{s}\}_{1
        \leq s \leq n^{i,\mathcal{L}}}$ and 
        $(\bar{\bu}^{i,\mathcal{L}},\bar{p}^{i,\mathcal{L}})$
        on $\Omega^{i,\ell}$ via MFEM
        
        \medskip

        \Return
        {
            $(\{\boldsymbol{\Phi}^{i,\mathcal{L}}_{s},\Psi^{i,\mathcal{L}}_{s}\},\,
            \bar{\bu}^{i,\mathcal{L}},\bar{p}^{i,\mathcal{L}})$
        }
    }
    \Else
    {
    	Define $i_1 = 2i-1$ and $i_2 = 2i$
        
        \medskip
        
        Decompose $\Omega^{i,\ell} = \Omega^{i_1,\ell+1} \cup \Omega^{i_2,\ell+1}$
        
        \medskip

        $\{\boldsymbol{\Phi}^{i_1,\ell+1}_{k},\,\Psi^{i_1,\ell+1}_{k}\}_{1 \leq k
        \leq n^{i_1,\ell+1}}$ and 
        $(\bar{\bu}^{i_1,\ell+1},\bar{p}^{i_1,\ell+1})
        \leftarrow$ RecMRCM($\Omega^{i_1,\ell+1}$) 

        \medskip

        $\{\boldsymbol{\Phi}^{i_2,\ell+1}_{k},\,\Psi^{i_2,\ell+1}_{k}\}_{1 \leq k \leq n^{i_2,\ell+1}}$ 
        and 
        $(\bar{\bu}^{i_2,\ell+1},\bar{p}^{i_2,\ell+1})$
        $\leftarrow$ RecMRCM($\Omega^{i_2,\ell+1}$) 

        \medskip

        Compute the coefficients $X^{i,\ell}_s$ and $\bar{X}^{i,\ell}$ by solving 
        \eqref{eq:linsys} on $\gamma^{i,\ell}$.

        \medskip

        \If{$\ell \neq 0$}
        {
            Compute 
            $\{\boldsymbol{\Phi}^{i,\ell}_{s},\Psi^{i,\ell}_{s}\}_{1 \leq s \leq n^{i,\ell}}$ 
            and 
            $(\bar{\bu}^{i,\ell},\bar{p}^{i,\ell})$
            on $\Omega^{i,\ell}$ with 
            \eqref{eq:new_bfs.1}-\eqref{eq:new_bfs.2}   
            
            \medskip

            \Return
            {
                $(\{\boldsymbol{\Phi}^{i,\ell}_{s},\Psi^{i,\ell}_{s}\},\,
                \bar{\bu}^{i,\ell},\bar{p}^{i,\ell})$
            }
        }   
        \Else
        {
            Compute $(\bu_h,p_h)$ on $\Omega$ with  
            \eqref{eq:new_bfs.1}-\eqref{eq:new_bfs.2}   

            \medskip

            \Return
            {
                $(\bu_h,p_h)$
            }
        }
    }
}
\caption{Recursive formulation for the MRCM}
\label{alg:recursive}
\end{algorithm}
\vspace{3mm}

\begin{figure}[h!]
	\centering
	\includegraphics[width=1.0\textwidth]{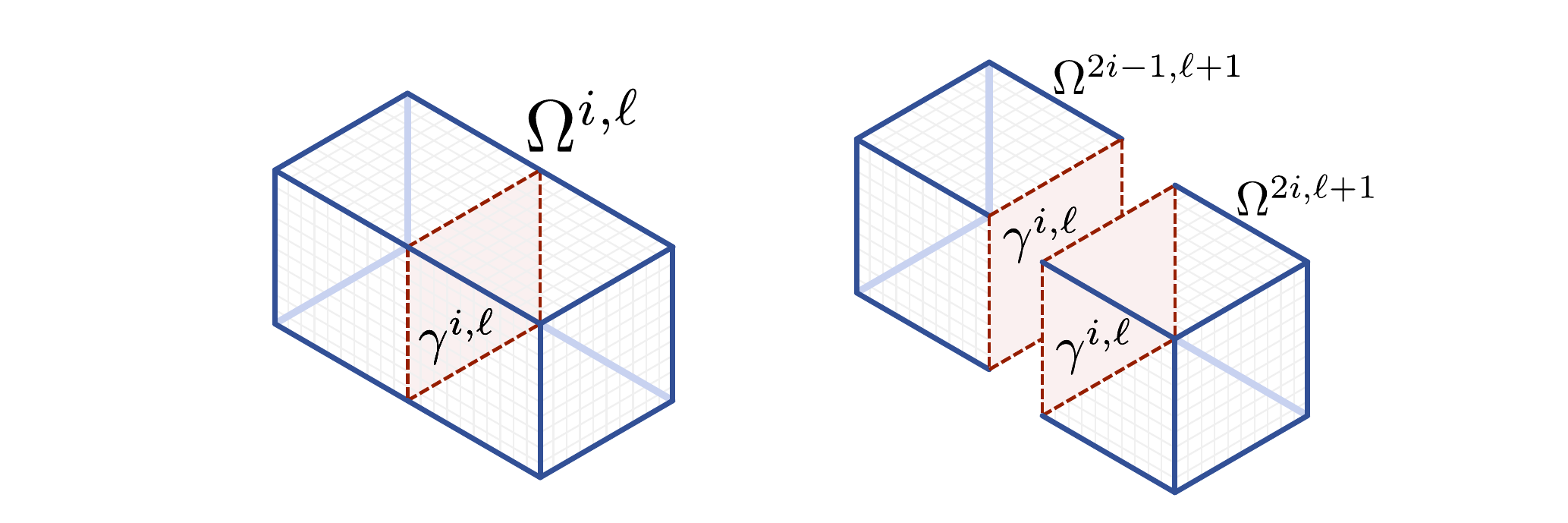}
	\caption{
        The representation of a subdomain $\Omega^{i,\ell}$ in level $\ell$, in which
        a step of the Recursive MRCM will be performed over $\gamma^{i,\ell}$.
        Local problems are defined in the new subdomains $\Omega^{2i-i,\ell+1}$
        and $\Omega^{2i,\ell+1}$.  
        }
	\label{fig:dd.4}
\end{figure}

We begin by defining
$(\bu^{i,\ell},p^{i,\ell})$ as the solution of the local problems
\eqref{eq11.2}-\eqref{eq14.2} restricted to $\Omega^{i,\ell}$.  
The solution is obtained by following the additive decomposition of the
MRCM, only now it is defined for each level: 
$\bu^{i,\ell} = \widehat{\bf u}^{i,\ell} + \bar{\bf u}^{i,\ell}$ and 
$p^{i,\ell} = \widehat{p}^{i,\ell} + \bar{p}^{i,\ell}$.
In each $\Omega^{i,\ell}$ we need to compute a set of associated MMBFs.
First, let us denote the set of MMBFs in $\Omega^{i,\ell}$ by
\begin{equation}\label{eq:mmbfs.rec.1}
(\boldsymbol{\Phi}^{i,\ell}_{s}\,,\Psi^{i,\ell}_{s}),
\quad
s = 1,\ldots,n^{i,\ell}, 
\quad
i=1,\ldots,m^{\ell}.
\end{equation}
At each level, the MMBFs are obtained by the solution of the local problems 
\eqref{eq11Hhatbasisa}-\eqref{eq12Hhatbasisb} on $\Omega^{i,\ell}$. 
Remember that a particular local solution,  
$(\bar{\bu}^{i,\ell},\bar{p}^{i,\ell})$, 
is also needed in order to complete the additive decomposition. 
The recursive MRCM algorithm can be described as follows: 
Consider $\Omega^{i,\ell}$, a generic subdomain of level $\ell$. 
We want to compute the MMBFs and the particular solution, 
$\{(\boldsymbol{\Phi}^{i,\ell}_{s}$,$\Psi^{i,\ell}_{s}),\,
(\bar{\bu}^{i,\ell},\bar{p}^{i,\ell})\}$,
$s = 1,\ldots,n^{i,\ell}$, associated with this subdomain. 
If $\ell \neq \mathcal{L}$, then we decompose $\Omega^{i,\ell}$
into $\Omega^{2i-i,\ell+1}$ and $\Omega^{2i,\ell+1}$,
as shown in Figure \ref{fig:dd.4}, 
and compute their associated MMBFs,
\begin{eqnarray}
\{(\boldsymbol{\Phi}^{2i-i,\ell+1}_{k}\,,\Psi^{2i-i,\ell+1}_{k}),\,
\bar{\bu}^{2i-i,\ell+1},\bar{p}^{2i-i,\ell+1}\}, 
\quad
&k = 1,\ldots,n^{2i-i,\ell+1} 
\label{eq:mmbfs_level.1} \\ 
\{(\boldsymbol{\Phi}^{2i,\ell+1}_{k}\,,\Psi^{2i,\ell+1}_{k}),\,
\bar{\bu}^{2i,\ell+1},\bar{p}^{2i,\ell+1}\}, 
\quad
&k = 1,\ldots,n^{2i,\ell+1},  
\label{eq:mmbfs_level.2}
\end{eqnarray}
in $\Omega^{2i-i,\ell+1}$ and $\Omega^{2i,\ell+1}$, respectively.

To compute \eqref{eq:mmbfs_level.1} and \eqref{eq:mmbfs_level.2} we apply the
MRCM again, decomposing each subdomain in level $\ell + 1$  
into two smaller subdomains and computing its MMBFs and particular solution
as above. This continues until $\ell = \mathcal{L}$, where
\eqref{eq:mmbfs_level.1} and \eqref{eq:mmbfs_level.2} are computed by MFEM.
It is important to notice that for $\ell = 0$ we do not compute MMBFs, but
the actual approximate solution $(\bu_h,p_h)$.
To complete the algorithm we need to compute, 
for each MMBF 
(and a particular solution) on $\Omega^{i,\ell}$ 
a set of coefficients $X^{i,\ell}$ (resp. $\bar{X}^{i,\ell}$, for the
particular local solution) in $\gamma^{i,\ell}$ 
by solving an interface linear system given by \eqref{eq:linsys} in the case  
of two subdomains. 
Then $(\boldsymbol{\Phi}^{i,\ell}_{s}\,,\Psi^{i,\ell}_{s},\,
\bar{\bu}^{i,\ell},\bar{p}^{i,\ell})$, $s = 1,\ldots,n^{i,\ell}$, 
are computed by a linear combination of 
\eqref{eq:mmbfs_level.1}-\eqref{eq:mmbfs_level.2} with its respective 
coefficients given by $X^{i,\ell}_s$ and $\bar{X}^{i,\ell}$. 

\begin{figure}[h!]
	\centering
	\includegraphics[width=1.0\textwidth]{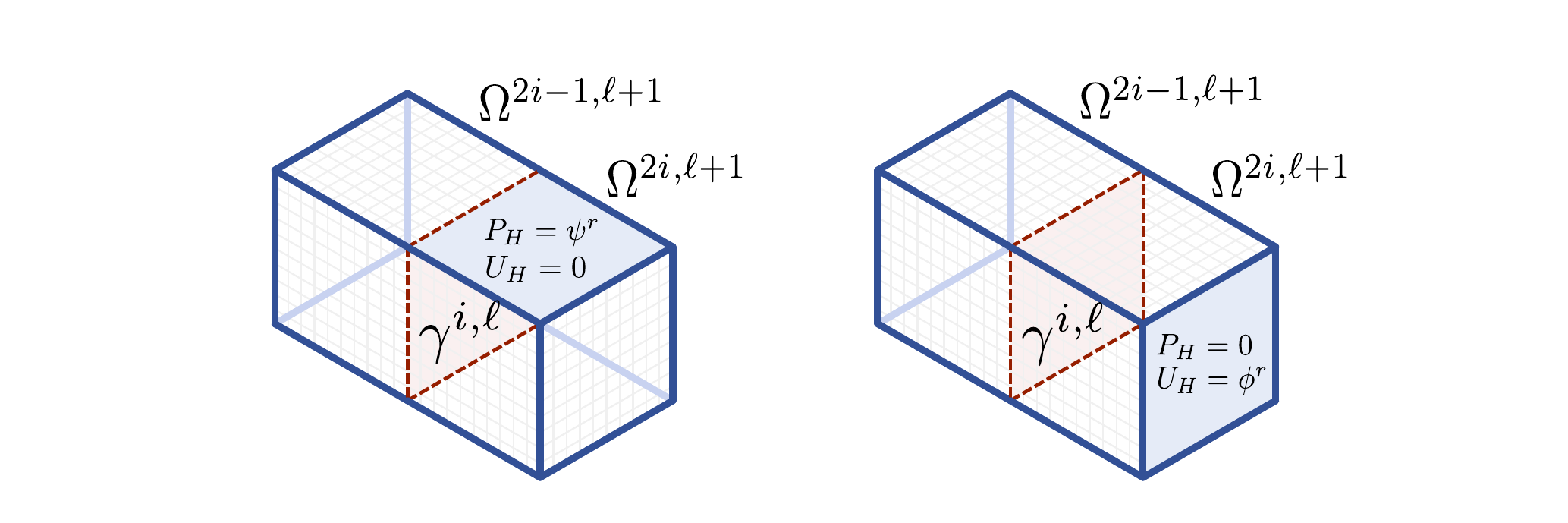}
	\caption{
        Representation of two boundary conditions of two different MMBFs
        in $\Omega^{i,\ell}$.
        This subdomain is composed of $\Omega^{2i-i,\ell+1}$ and $\Omega^{2i,\ell+1}$.
        The MMBFs of $\Omega^{i,\ell}$ are obtained by alternating the values of $U_H$ 
        and $P_H$ on the external interfaces boundaries $\Gamma^{i,\ell}$, 
        for their finite element basis $\phi^j$ and $\psi^j$, as in
        \eqref{eq11Hhatbasisa}-\eqref{eq12Hhatbasisb}.
        On the left figure, the value of $P_H = \psi^r$ is the contribution of that 
        particular coarse interface.
        On the right figure, the value of $U_H = \phi^r$ is the contribution of that 
        particular coarse interface.
        The MMBFs in $\Omega^{i,\ell}$ are obtained by the linear combinations 
        \eqref{eq:new_bfs.1}-\eqref{eq:new_bfs.2}.
	}
	\label{fig:new_mmbf}
\end{figure}

Let us now explain how to compute the MMBFs from 
the linear combination of MMBFs of subsequent levels. 
Suppose we already computed the MMBFs of $\Omega^{2i-i,\ell+1}$ 
and $\Omega^{2i,\ell+1}$ and its associated 
coefficients $X^{i,\ell},\bar{X}^{i,\ell}$ on level $\ell+1$. 
Then, each MMBFs on $\Omega^{i,\ell}$ is computed by 
\begin{eqnarray}
\boldsymbol{\Phi}^{i,\ell}_{s} &=  
\sum\limits_{k} 
X^{i,\ell}_{k,\,s} \boldsymbol{\Phi}^{2i-1,\,\ell+1}_{k}
+ \bar{\bu}^{2i-1,\,\ell+1} 
+ \, \sum\limits_{k} 
X^{i,\ell}_{k,\,s} \boldsymbol{\Phi}^{2i,\,\ell+1}_{k}
+ \bar{\bu}^{\ell+1\,,2i}  
+  \phi^r\, \boldsymbol{\Phi}^{2i,\,\ell+1}_{r}, 
\label{eq:new_bfs.1} \\
[3mm]
\Psi^{i,\ell}_{s} &=  
\sum\limits_{k} 
X^{i,\ell}_{k,\,s} \Psi^{2i-1,\,\ell+1}_{k}
+ \bar{p}^{2i-1,\,\ell+1} 
+ \, \sum\limits_{k} 
X^{i,\ell}_{k,\,s} \Psi^{2i,\,\ell+1}_{k}
+ \bar{p}^{2i,\,\ell+1}
+  \phi^r\, \Psi^{2i,\,\ell+1}_{r}, 
\label{eq:new_bfs.2}
\end{eqnarray}
where $s \in \{1,\ldots,n^{i,\ell}\}$ and $k \in
\{1,\ldots,\xi^{i,\ell}\}$; $\xi^{i,\ell}$ is the number of
interface degrees of freedom on $\gamma^{i,\ell}$.
The last terms in \eqref{eq:new_bfs.1} and \eqref{eq:new_bfs.2} are related to
the MMBF ($\boldsymbol{\Phi}^{2i,\ell+1}_{r},\,\Psi^{2i,\ell+1}_{r}$) on
$\Omega^{2i,\ell+1}$ 
that accounts for the contribution of the boundary value $\phi^r$ in $\Gamma^{i,\ell}$, 
as illustrated in the right figure in Figure \ref{fig:new_mmbf}.
The particular solution can also be written as a linear combination similar to 
\eqref{eq:new_bfs.1}-\eqref{eq:new_bfs.2},
\begin{eqnarray}
\bar{\bu}^{i,\ell} &=  
\sum\limits_{k} 
\bar{X}^{i,\ell}_{k} \boldsymbol{\Phi}^{2i-1,\,\ell+1}_{k}
+ \bar{\bu}^{2i-1,\,\ell+1} 
+ \, \sum\limits_{k} 
\bar{X}^{i,\ell}_{k} \boldsymbol{\Phi}^{2i,\,\ell+1}_{k}
+ \bar{\bu}^{2i,\,\ell+1},  
\label{eq:particular.1} \\
[3mm]
\bar{p}^{i,\ell} &=  
\sum\limits_{k} 
\bar{X}^{i,\ell}_{k} \Psi^{2i-1,\,\ell+1}_{k}
+ \bar{p}^{2i-1,\,\ell+1} 
+ \, \sum\limits_{k} 
\bar{X}^{i,\ell}_{k} \Psi^{2i,\,\ell+1}_{k}
+ \bar{p}^{2i,\,\ell+1},
\label{eq:particular.2}
\end{eqnarray}
where $k \in \{1,\ldots,\xi^{i,\ell}\}$.

{\bf Remark:}
This recursive construction of the approximate solution of 
\eqref{eq:ellip.1}-\eqref{eq:ellip.3} by the MRCM allows us to decompose
the global interface linear system \eqref{eq:linsys} into a set of small
and localized interface linear systems on $\gamma^{i,\ell}$ for all subdomains $i$
and all levels $\ell$.
Each set of local linear systems on $\gamma^{i,\ell}$ has size $\xi^{i,\ell}
\times \xi^{i,\ell}$. The
linear systems are independent of each other and can be solved
simultaneously.
The matrix and right hand side of the local interface systems are constructed
using \eqref{eq:linsys} restricted to $\gamma^{i,\ell}$. 
One important aspect of the algorithm
described here is that it can keep track of the coefficients of the linear
combinations that are used to express each MMBF of each level as a linear
combination of the MMBfs of the previous levels, in a way that we do not need
to store all the values of the MMBFs in \textit{coarser} levels.
As we proceed to \textit{coarser} levels,
those MMBFs can be expressed as linear combinations of the MMBFs
associated with the \textit{finest} decomposition ($\ell = \mathcal{L}$). 
This way, for any given level, we are able to express each multiscale 
basis function as a linear combination of the \textit{finest} level MMBFs.

% MuMM connection:
%\input{section/section_4}
\section{Parallel Implementation and connection to the Multiscale Mixed Method}\label{sec:3.4}

\begin{figure}[h!]
	\centering
    \includegraphics[width=0.3\textwidth]{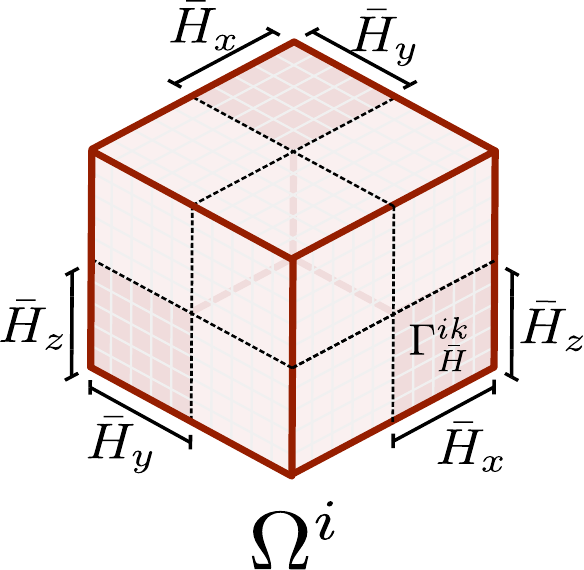}	
	\caption{
		Representation of $\bar H$ scale. 
	}
	\label{fig:bfs}
\end{figure}

The flexibility in the choice of interface spaces for pressure and normal
fluxes provided by the MRCM framework comes with a cost. Even if piecewise
constant spaces are selected for both variables in a three-dimensional (resp.
two-dimensional) subdomain the minimum number of MMBFs that need to be computed
in each subdomain
is 12 (resp. 8). The problems we intend to solve using the recursive framework
will involve up to billions of cells. Thus, we wish to perform simulations with
methods that are as inexpensive as possible from the computational point of
view. In this context we will implement the recursive procedure in a particular
case of MRCM: the Multiscale Mixed Method (MuMM).  In implementing the MuMM one
need only a set of six (resp. four) MMBFs, in three (resp. two) dimensions,
thus reducing the computational cost of the implementation.
An important feature of the MuMM is the introduction of an intermediate
coarse scale of size $\bar{H}$, such that $h \leq \bar{H} \leq H$, where we define the 
interface space $F_{\bar H} \subset F_h(\mathcal{E}_h)$. 
This space is taken to be piecewise constant in the $\bar H$ 
scale, see Figure \ref{fig:bfs}.

Let us define $\Gamma^{ik}_{\bar H}$ as one element of the partition of $\Gamma$, with
size $\bar H$, adjacent to subdomains $\Omega^i$ and $\Omega^k$, such that 
$H/\bar{H}$ and $\bar{H}/h$ are both integer numbers. 
This partition can also be performed
independently for each direction, with minor modifications.
The introduction of an intermediate scale does not changes the recursive formulation 
construction, it only adds a flexibility in the number of MMBFs
and in the size of the interface linear system to be computed. 
Since the continuity equations in $\Gamma$ 
are defined in the coarse scale, flux conservation is 
only satisfied in this scale.
Downscaling (or smoothing) techniques should be used to recover flux conservation
on fine scale \cite{Pereira2014,APGK19,guiraldello2018multiscale,Rafael2020}. 

%The recursive formulation was implemented in C version 6.5, C++ version 6.5 
%and openMPI version 2.0.4.2.
The recursive formulation was implemented in C, C++ and openMPI.
To compute the MMBFs and particular for each subdomain of level $\ell = \mathcal{L}$, 
we use a Mixed Finite Element discretization with lowest index Raviart-Thomas  
spaces \cite{RT91} to construct a linear system for the pressure variable. 
The solution was obtained by means of a conjugate gradient with
an algebraic multigrid preconditioner C++ solver, with a tolerance of $10^{-8}$ 
\cite{AMGtoolbox}.  
The interface linear systems were solved by a simple,
in-house implemented LU solver 
since its matrix can be quite small (depending on the choice of the size of the
$\bar H$ scale) and are efficiently computed by such solver. 
The recursive formulation is implemented considering a decomposition of the
domain such that 
each direction is decomposed in a power of two. 
This simplifies the implementation of the message passing between subdomains.

\begin{figure}[h!]
	\centering
    \includegraphics[width=1.\textwidth]{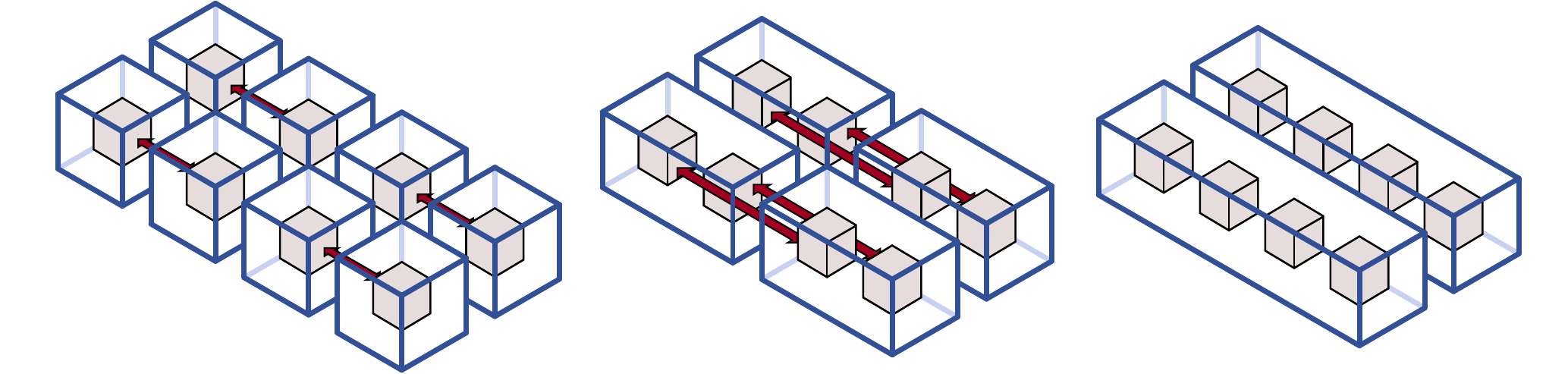}
    \caption{Representation of the message passing pattern between cores (red arrows) in the
    union operations, as seen in Section 3. In the first level
    the cores (represented by the grey cubes) and subdomain (represented by
    shallow blue cubes) meshes are the same. The communication is done with
    their direct adjacent subdomain. 
    However, in the coarse levels the cores that compose a new subdomain
    communicate with the cores that are in the same ``position'' inside 
    the new adjacent subdomain.  
    }
	\label{fig:cluster.3}
\end{figure}

The exchange of information between subdomains is done by
keeping the same number of message passing steps constant at each level. 
This is achieved by a one-to-one message passing between the processing cores that
compose a subdomain at a given level. This is illustrated in Figure
\ref{fig:cluster.3}. 
All our experiments were done on the Santos Dumont cluster 
located at the National Laboratory for Scientific Computing
(LNCC) in Petr\'opolis, RJ, Brazil, from several million to $2$ billion cells
on a dual-socket Intel Xeon E5-2695-v2, 2,4GHZ, 64GB DDR3 RAM.
    %

% Numerical experiments:
%\input{section/section_5}
\section{Numerical Experiments: Setup, Results and Discussions}\label{sec:3.5}

\subsection{Setup for the numerical experiments}

In this section, we present numerical experiments to evaluate the computational
efficiency and the accuracy of solution of our three-dimensional parallel
implementation of the recursive formulation for very large problems, up to $2$
billion cells.   

We consider the pressure-velocity problem \eqref{eq:ellip.1}-\eqref{eq:ellip.3}, 
for a physical domain $[0, L_x] \times [0, L_y] \times [0, L_z]$
and isotropic absolute permeability tensors.
Our implementation is based on the MuMM \cite{Pereira2014}, 
where the interface spaces are piecewise constant functions, as explained in
Section \ref{sec:3.4}. 
For every $\Gamma^{ik}_{\bar H}$, the Robin parameter $\beta^i$ and $\beta^k$
are chosen to be constant both defines as
\begin{equation}\label{eq:robin.beta}
\beta^i = \beta^k = \frac{\alpha\,\bar{H}}{\bar{K}_{eff}}, 
\end{equation}
where $\alpha$ is a dimensionless parameter
\cite{guiraldello2018multiscale,Pereira2014,Guiraldello2019} and
$\bar{K}_{eff}$ is the average of all harmonic means of the adjacent permeabilities 
in the cells that compose the $\bar H$ scale, i.e.,
\begin{equation}\label{eq:permeff}
\bar{K}_{eff} = \frac{1}{N} \sum\limits_{e \in \Gamma^{ik}_{\bar H}}
~\frac{2 K^i_e K^k_e}{K^i_e + K^k_e},
\end{equation}
the sum is on all $e$ cells that compose the $\Gamma^{ik}_{\bar H}$.

As discussed before, the magnitude of $\alpha$ controls the coupling between
the subdomains, as explained in \cite{guiraldello2018multiscale}.
The choice of large $\alpha$ values gives higher priority to flux continuity
over pressure continuity of the final solution. For the scalability studies
we choose a constant value $\alpha = 10^{3}$ and for the accuracy results 
we choose $\alpha = 10^{6}$.
The number of operations of the recursive algorithm remains the same 
and therefore we do not expect its scalability to be affected by the $\alpha$ value
(for more details on how the magnitude of $\alpha$ affects the solution, see
\cite{Guiraldello2019}). 

The computational efficiency is evaluated in two cases, namely: $i)$
homogeneous permeability field and $ii)$ high contrast heterogeneous
permeability field. For both cases we perform scalability studies where we
assess the behavior of the numerical method relative to its computational time
against an increasing number of cores.  The first scalability test is the
strong scaling, where the total number of discretization elements and problem
size is fixed while increasing the total number of processing cores.  
The second scalability test is the weak scaling, 
where the total size of the problem is increased, while increasing the number
of processing cores. We keep the size of the local linear systems in each
subdomain constant, while increasing the overall problem size and processing
cores, therefore the expected simulation time should remain constant throughout
the weak scaling tests. 
The boundary conditions are given by $p(0,y,z) = 1$ and $p(L_x,y,z) = 0$, combined
with no-flow conditions on the other boundaries. For the weak scaling case,
boundary conditions are updated in each case, in order to keep the same overall
flux, making sure the Darcy problem stays the same, at least for the
homogeneous permeabilities.

In the recursive algorithm we need to establish the mesh for the 
finest domain decomposition,
so that we associate each subdomain in level $\ell = \mathcal{L}$ with a unique
core. In all our experiments, we define the finest decomposition with no divisions on
the $z$-direction. 
The implementation considers domain decompositions where
each direction is decomposed in a power of two. 
Tables \ref{tab:experiments1} and \ref{tab:experiments2} organizes the scaling experiments, 
showing the subdomain divisions and number of cells 
for the strong and weak scaling studies.
% Please add the following required packages to your document preamble:
% \usepackage{multirow}
\begin{table}[h!]
\centering
\caption{Setup for the scaling experiments with up to $134$ million cells (strong
scaling) and $268$ million cells (weak scaling).}
\sizef
\begin{tabular}{rc|c|rc|c}
\hline
\hline
\multicolumn{3}{c|}{\multirow{2}{*}{Strong Scaling}} & \multicolumn{3}{c}{\multirow{2}{*}{Weak Scaling}} \\
\multicolumn{3}{c|}{} & \multicolumn{3}{c}{} \\ \hline \hline
\multicolumn{3}{c|}{\multirow{2}{*}{$\sim$ 134 million in $\Omega$}} 
& \multicolumn{3}{c}{\multirow{2}{*}{$\sim$ 262 thousand in each $\Omega^i$}} \\
\multicolumn{3}{c|}{} & \multicolumn{3}{c}{} \\ \hline
\multicolumn{2}{c|}{Cores (subdomains)} & Total cells in $\Omega^i$ & \multicolumn{2}{c|}{Cores (subdomains)} & Total cells in $\Omega$ \\ \hline
32 & (4 $\times$ 8 $\times$ 1) & 128 $\times$ 64 $\times$ 512 & 32 & (4 $\times$ 8 $\times$ 1)        & $8.39 \times 10^6$ \\
64 & (8 $\times$ 8 $\times$ 1) & 64 $\times$ 64 $\times$ 512 & 64 & (8 $\times$ 8 $\times$ 1)         & $1.68 \times 10^7$ \\
128 & (8 $\times$ 16 $\times$ 1) & 64 $\times$ 32 $\times$ 512 & 128 & (8 $\times$ 16 $\times$ 1)     & $3.36 \times 10^7$ \\
256 & (16 $\times$ 16 $\times$ 1) & 32 $\times$ 32 $\times$ 512 & 256 & (16 $\times$ 16 $\times$ 1)   & $6.71 \times 10^7$ \\
512 & (16 $\times$ 32 $\times$ 1) & 32 $\times$ 16 $\times$ 512 & 512 & (16 $\times$ 32 $\times$ 1)    & $1.34 \times 10^8$ \\
1024 & (32 $\times$ 32 $\times$ 1) & 16 $\times$ 16 $\times$ 512 & 1024 & (32 $\times$ 32 $\times$ 1) & $2.68 \times 10^8$ \\ \hline
\end{tabular}
\label{tab:experiments1}
\end{table}
\begin{table}[h!]
\centering
\caption{Setup for the scaling experiments with up to $1$ billion cells (strong
scaling) and $2$ billion cells (weak scaling).}
\sizef
\begin{tabular}{rc|c|rc|c}
\hline
\hline
\multicolumn{3}{c|}{\multirow{2}{*}{$\sim$ 1 billion in $\Omega$}} &
\multicolumn{3}{c}{\multirow{2}{*}{$\sim$ 2 million in each $\Omega^i$}} \\
\multicolumn{3}{c|}{} & \multicolumn{3}{c}{} \\ \hline
\multicolumn{2}{c|}{Cores (subdomains)} & Total cells in $\Omega^i$ & \multicolumn{2}{c|}{Cores (subdomains)} & Total cells in $\Omega$ \\ \hline
\multirow{2}{*}{256} & \multirow{2}{*}{(16 $\times$ 16 $\times$ 1)} & \multirow{2}{*}{64 $\times$ 64 $\times$ 1024} & 32 & (4 $\times$ 8 $\times$ 1) & $6.71 \times 10^7$ \\
 &  &  & 64 & (8 $\times$ 8 $\times$ 1) & $1.34 \times 10^8$ \\
\multirow{2}{*}{512} & \multirow{2}{*}{(16 $\times$ 32 $\times$ 1)} & \multirow{2}{*}{64 $\times$ 32 $\times$ 1024} & 128 & (8 $\times$ 16 $\times$ 1) & $2.68 \times 10^8$ \\
 &  &  & 256 & (16 $\times$ 16 $\times$ 1) & $5.37 \times 10^8$ \\
\multirow{2}{*}{1024} & \multirow{2}{*}{(32 $\times$ 32 $\times$ 1)} &
\multirow{2}{*}{32 $\times$ 32 $\times$ 1024} & 512 & (16 $\times$ 32 $\times$ 1) & $1.07\times 10^9$ \\
 &  &  & 1024 & (32 $\times$ 32 $\times$ 1) & $2.15 \times 10^9$ \\ \hline
\end{tabular}
\label{tab:experiments2}
\end{table}

Next, we need to define the size of the coarse $\bar H$ partition.
For the experiments we use two sets of coarse scale in each 
$\Omega^{i,\mathcal{L}}$: 
$\bar{H}_x = H_x$, $\bar{H}_y = H_y$, and 
$\bar{H}_x = H_x/4$, $\bar{H}_y = H_y/4$;  
for the coarse scale in the $z$-direction, we fixed $\bar{H}_z = H_z$.
The number of $\bar{H}$ partitions on $\Gamma^{i,\mathcal{L}}$ 
is the number of MMBFs to be directly computed on the finest decomposition
by HMFEM, and it defines the number
of MMBFs on all levels through \eqref{eq:new_bfs.1}-\eqref{eq:new_bfs.2}.
Table \ref{tab:ntmmbfs} shows, for each coarse scale partition chosen, 
the total number of MMBFs and particular solution to be computed with HMFEM 
for all subdomains in the last level.
The direct computation of local problems by HMFEM is the most expensive part of
the algorithm, as we will see in the experiments below. 
The last column shows the increase percentage in the total number of local
problems to be computed. 
\begin{table}[h!]
\centering
\caption{Shows the total number (globally) of MMBFs to be computed in level
$\ell = \mathcal{L}$ for the cases where we have $\bar{H} = H$ and $\bar{H} =
H/4$.}
\sizef
\begin{tabular}{c|c|c|c}
\hline\hline
Subdomains & \begin{tabular}[c]{@{}c@{}}Number of MMBFS for \\ $\overline{H}=H$\end{tabular} & \begin{tabular}[c]{@{}c@{}}Number of MMBFS for \\ 
$\overline{H}=H/4$\end{tabular} & $\%$ increase \\ \hline \hline
4 $\times$ 8 $\times$ 1   &  136 &   448 & 330 $\%$  \\
8 $\times$ 8 $\times$ 1   &  288 &   960 & 333 $\%$  \\
8 $\times$ 16 $\times$ 1  &  592 &  1984 & 335 $\%$  \\
16 $\times$ 16 $\times$ 1 & 1216 &  4096 & 337 $\%$  \\
16 $\times$ 32 $\times$ 1 & 2464 &  8320 & 338 $\%$  \\
32 $\times$ 32 $\times$ 1 & 4992 & 16896 & 338 $\%$  \\ \hline
\end{tabular}
\label{tab:ntmmbfs}
\end{table}

%Values not in percentage
%& 3.2941 \\
%& 3.3333 \\
%& 3.3514 \\
%& 3.3684 \\
%& 3.3766 \\
%& 3.3846 \\

Finally, in the accuracy experiments we show that the accuracy of the
approximated flux does not deteriorates as we increase the number of cores for
the strong and weak scaling studies. No downscaling strategy was used so
that pressure and normal fluxes may be discontinuous at the fine grid across
the skeleton of the decomposition.
As we are dealing with very large problems, we restricted our simulations to a
maximum of $4$ million cells per subdomain due to memory and computational
restrictions.  

\subsection{Homogeneous scalability study}

For the experiments in this section, 
we consider an isotropic homogeneous absolute permeability field
given by $K(\bx) = 1$. 

\subsubsection{Strong scaling}

\begin{figure}[h!]
	\centering
	\setlength{\abovecaptionskip}{-35pt}
	\subfloat[$134$ million cells.]{ \label{fig:ss.ho.3d.a}\includegraphics[width=0.47\textwidth]{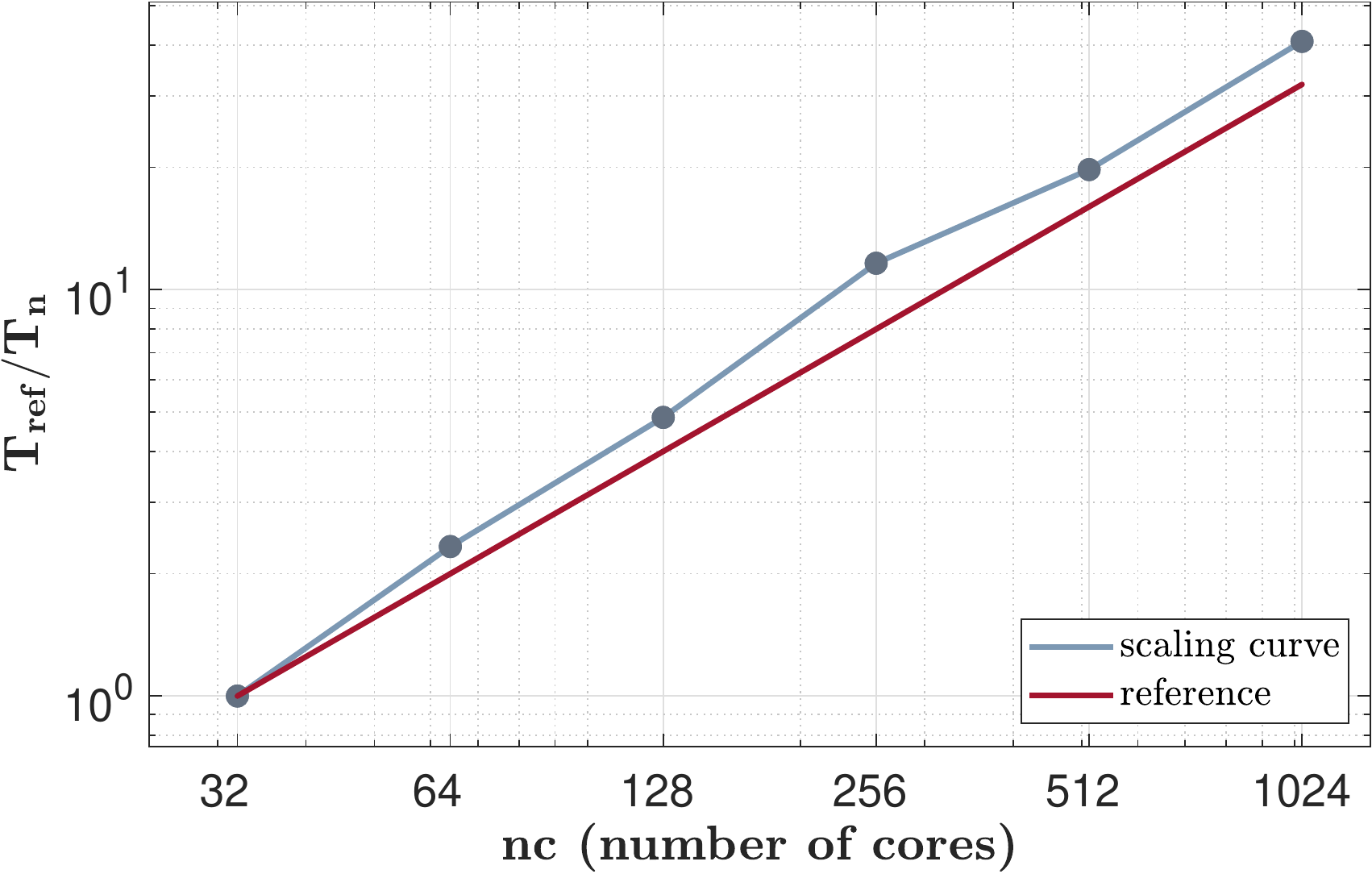}}
	\qquad
	\subfloat[$1$ billion cells.]{ \label{fig:ss.ho.3d.b}\includegraphics[width=0.46\textwidth]{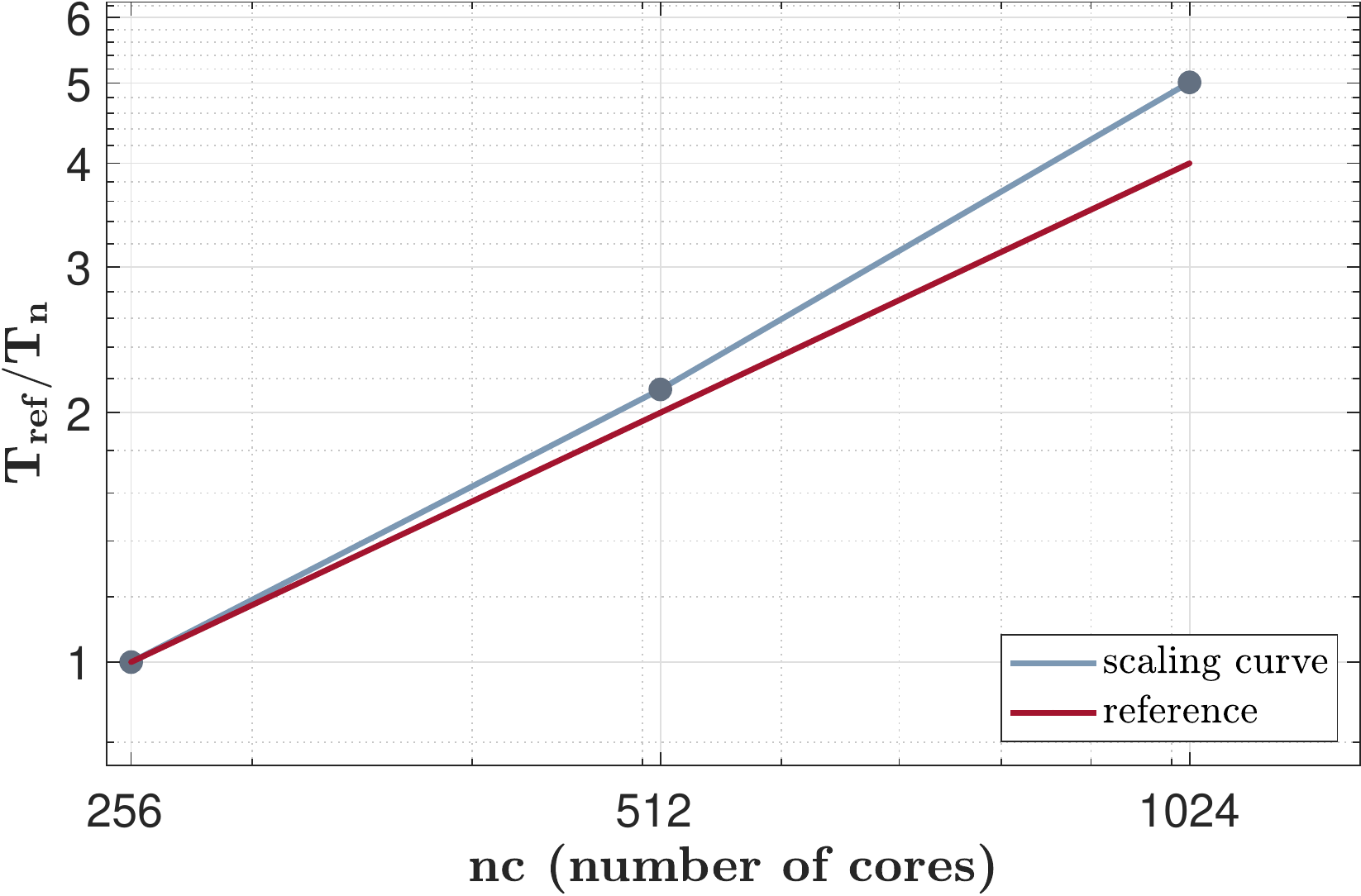}}
	\vspace{15mm}
	\caption{
		%\footnotesize
        Strong scaling curves with homogeneous permeability and $\bar{H} = H$ (see Table \ref{table:ss.ho.3d}). 
	}
	\label{fig:ss.ho.3d}
\end{figure}

\begin{figure}[h!]
	\centering
	\setlength{\abovecaptionskip}{-35pt}
	\subfloat[$134$ million cells.]{\label{fig:ss.ho4.3d.a}\includegraphics[width=0.47\textwidth]{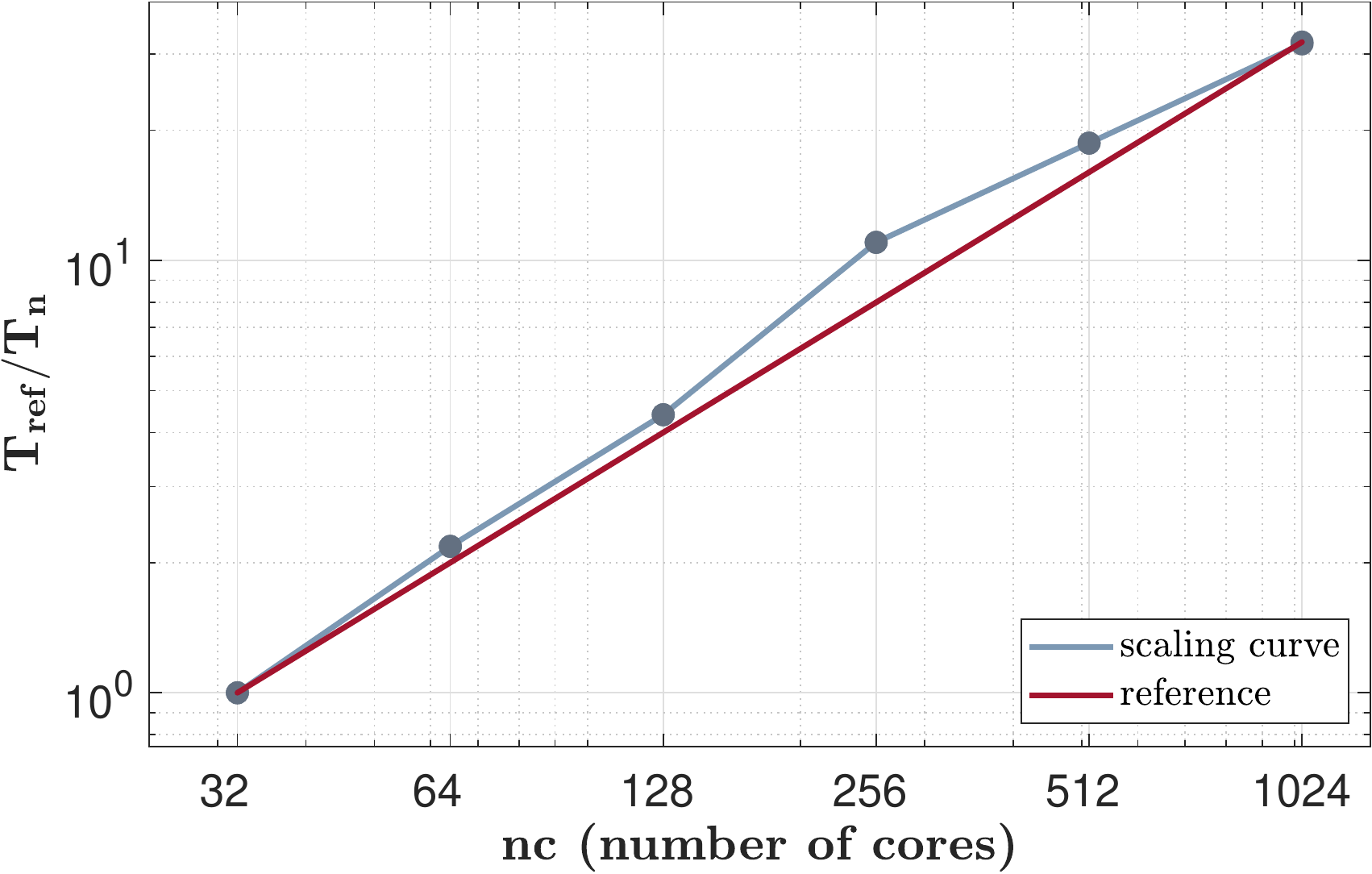}}\qquad
	\subfloat[$1$ billion cells.]{\label{fig:ss.ho4.3d.b}\includegraphics[width=0.46\textwidth]{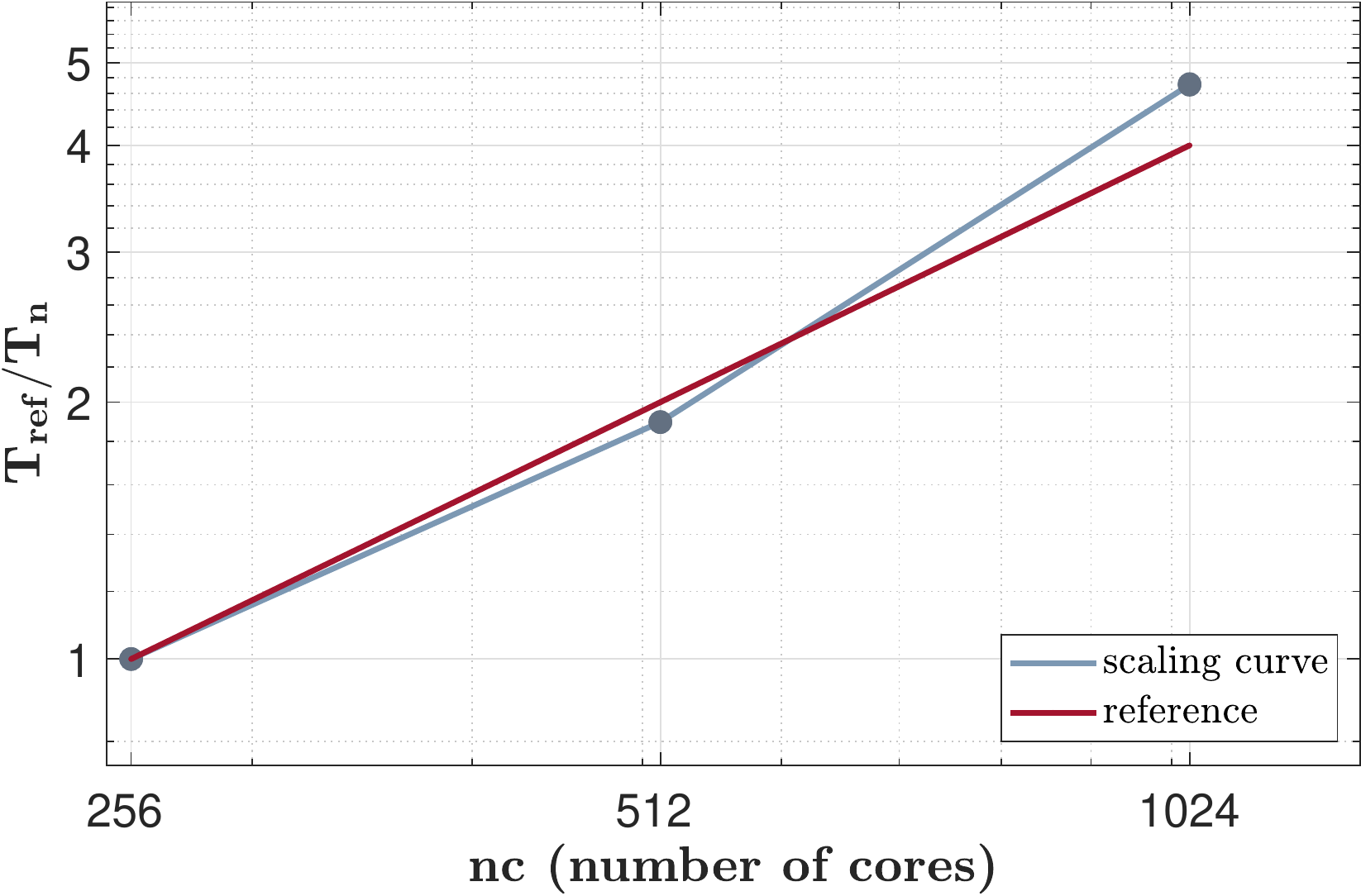}}
	\vspace{15mm}
	\caption{
		%\footnotesize
        Strong scaling curves with homogeneous permeability and $\bar{H} = H/4$ (see Table \ref{table:ss.ho4.3d}). 
	} 
	\label{fig:ss.ho4.3d}
\end{figure}

Figures \ref{fig:ss.ho.3d} and \ref{fig:ss.ho4.3d} 
present the scalability curves of time ratio versus number of cores.
Under ideal conditions, with a fully paralellizable method,
we expect the computational time to be 
reduced by half if we double the number of processing cores, since
the computational power was doubled. 
%However, according to the Amadahl's Law, the speedup gain of any algorithm is
%limited by the part of that algorithm that is not parallelizable \cite{amdahl67}.
The red curve represents the ideal scaling curve,
\begin{equation}\label{eq:stg_scal_curve}
\frac{T_{ref}}{T_{n}} = \frac{n_c}{n_{c_{ref}}},
\end{equation}
where $T_n$ is the total time of computation and $n_c$ the number of cores used;
while $T_{ref}$ is the reference processing time spent to compute the solution
using $n_{c_{ref}}$ cores.  The blue curves represent our data.

In Figures \ref{fig:ss.ho.3d.a} and \ref{fig:ss.ho.3d.b} we present the scaling
curves for the three-dimensional problem with $512 \times 512 \times 512$
(approximately $134$ million) cells and $1024 \times 1024 \times 1024$
(approximately $1$ billion) cells, respectively, with computational times
reported in Table \ref{table:ss.ho.3d} for the $\bar{H} = H$ case.  
The same experiments
are reported in Figures \ref{fig:ss.ho4.3d.a} and \ref{fig:ss.ho4.3d.b} as well
as in Table \ref{table:ss.ho4.3d}, for the $\bar{H} = H/4$ case.
\begin{table}[h!]
\centering
\caption{
Strong scaling times for homogeneous problem
with $134$ million cells (top table, see Figure \ref{fig:ss.ho.3d.a}) and 
with $1$ billion cells (bottom table, see Figure \ref{fig:ss.ho.3d.b}). 
For these problems we considered 
$\bar{H} = H$. 
}
\sizefonte
%\begin{tabular}{\textwidth}{c|c|c|c|c}
\begin{tabular}{c|c|c|c|c|c}
\hline \hline
\multicolumn{6}{c}{\multirow{2}{*}{\textbf{Strong Scaling (homogeneous permeability -  $\bar{H} = H$)}}} \\ 
\multicolumn{6}{c}{} \\ 
\hline \hline
\multicolumn{6}{c}{\multirow{2}{*}{\textbf{134 million cells}}} \\ 
\multicolumn{6}{c}{} \\ 
\hline
Cores  &  MMBFs Time (s) & INTRF Time (s) & MPI Time (s) & Total Time (s) & $\%$ decrease (Total) \\ \hline
32     & 94.66 & 0.0071 & 0.0001 & 128.84 &          \\ 
64     & 41.43 & 0.0072 & 0.0002 & 54.88  & 57.41 \\
128    & 20.81 & 0.0048 & 0.0003 & 26.15  & 52.35  \\ 
256    & 8.31  & 0.0063 & 0.0005 & 10.75  & 58.87  \\ 
512    & 4.73  & 0.0082 & 0.0007 & 5.98   & 44.40  \\ 
1024   & 1.88  & 0.0123 & 0.0024 & 2.73   & 54.30 \\ \hline
\end{tabular}
%\begin{tabular}{\textwidth}{c|c|c|c|c}
\begin{tabular}{c|c|c|c|c|c}
\multicolumn{6}{c}{\multirow{2}{*}{\textbf{1 billion cells}}} \\ 
\multicolumn{6}{c}{} \\ 
\hline
Cores &  MMBFs Time (s) & INTRF Time (s) & MPI Time (s) & Total Time (s) & $\%$ decrease (Total) \\ \hline
256   & 87.107  & 0.1009 & 0.0062 & 116.22 &         \\ 
512   & 46.195  & 0.0108 & 0.0007 & 54.094 &  53.46   \\ 
1024  & 17.278  & 0.0151 & 0.0024 & 22.858 &  57.74   \\ \hline
\end{tabular}
\label{table:ss.ho.3d}
\end{table}

\begin{table}[h!]
\centering
\caption{
Strong scaling times for homogeneous problem 
with $134$ million cells (top table, see Figure \ref{fig:ss.ho4.3d.a}) and 
with $1$ billion cells (bottom table, see Figure \ref{fig:ss.ho4.3d.b}). 
For these problems we considered 
$\bar{H} = H/4$. 
}
\sizefonte
%\begin{tabularx}{\textwidth}{c|c|Y|Y|Y}
\begin{tabular}{c|c|c|c|c|c}
\hline \hline
\multicolumn{6}{c}{\multirow{2}{*}{\textbf{Strong Scaling (homogeneous permeability - $\bar{H} = H/4$)}}} \\ 
\multicolumn{6}{c}{} \\ 
\hline \hline
\multicolumn{6}{c}{\multirow{2}{*}{\textbf{134 million cells}}} \\ 
\multicolumn{6}{c}{} \\ 
\hline
Cores &  MMBFs Time (s) & INTRF Time (s) & MPI Time (s) & Total Time (s) & $\%$ decrease (Total) \\ \hline
32    & 303.72  & 0.0226 & 0.0013 & 396.36   &     \\ 
64    & 140.92  & 0.0240 & 0.0024 & 181.40   & 54.23      \\ 
128   &  72.60  & 0.0553 & 0.0037 &  89.74   & 50.53     \\ 
256   &  29.26  & 0.1596 & 0.0124 &  35.66   & 60.27     \\ 
512   &  16.46  & 0.9026 & 0.0356 &  20.66   & 42.06     \\ 
1024  &   6.71  & 3.9744 & 0.1262 &  12.12   & 41.34     \\ \hline
\end{tabular}
%\begin{tabularx}{\textwidth}{c|c|Y|Y|Y}
\begin{tabular}{c|c|c|c|c|c}
\multicolumn{6}{c}{\multirow{2}{*}{\textbf{1 billion cells}}} \\ 
\multicolumn{6}{c}{} \\ 
\hline
Cores &  MMBFs Time (s) &  INTRF Time (s) & MPI Time (s) & Total Time (s) & $\%$ decrease (Total) \\ \hline
256   & 297.61  & 0.6224 & 0.0410 & 362.87   &      \\
512   & 166.59  & 0.8838 & 0.0337 & 191.05   & 47.35    \\ 
1024  &  61.10  & 3.9706 & 0.1239 &  76.67   &  59.87   \\ \hline
\end{tabular}
\label{table:ss.ho4.3d}
\end{table}

We can see from these figures that our simulations
are above the optimal (red) curve up to $1024$ cores, 
showing outstanding parallel performance.
This behavior can be observed in Tables
\ref{table:ss.ho.3d} and \ref{table:ss.ho4.3d}, which show
the times for the total run and for the MMBFs computation 
on level $\ell = \mathcal{L}$. 
From these tables, one can see that the most expensive part of the overall
computation of the solution is in the construction of the MMBF's, as compared
to the total time. The remaining time includes the solution of the interface
problems in all levels of the recursive algorithm as well as the time spent
with exchange of information between subdomains. 
The latter is around three orders of magnitude smaller as compared to the
processing times.

%\newpage

\subsubsection{Weak scaling} 

As mentioned before, in an ideal problem ($100\%$ parallelizable) we expect a
constant overall processing time while increasing problem size, since the
degrees of freedom for each local problem are fixed. Scalability curves are
reported in Figures 9 and 10, where we can see overall processing time versus
number of cores. The reference curve (red curve) is an average of the total
times obtained by our simulations, reported by the blue curves.

\begin{figure}[h!]
	\centering
	\setlength{\abovecaptionskip}{-35pt}
	\subfloat[$262$ thousand cells per subdomain.]{ \label{fig:ws.ho.3d.a}\includegraphics[width=0.46\textwidth]{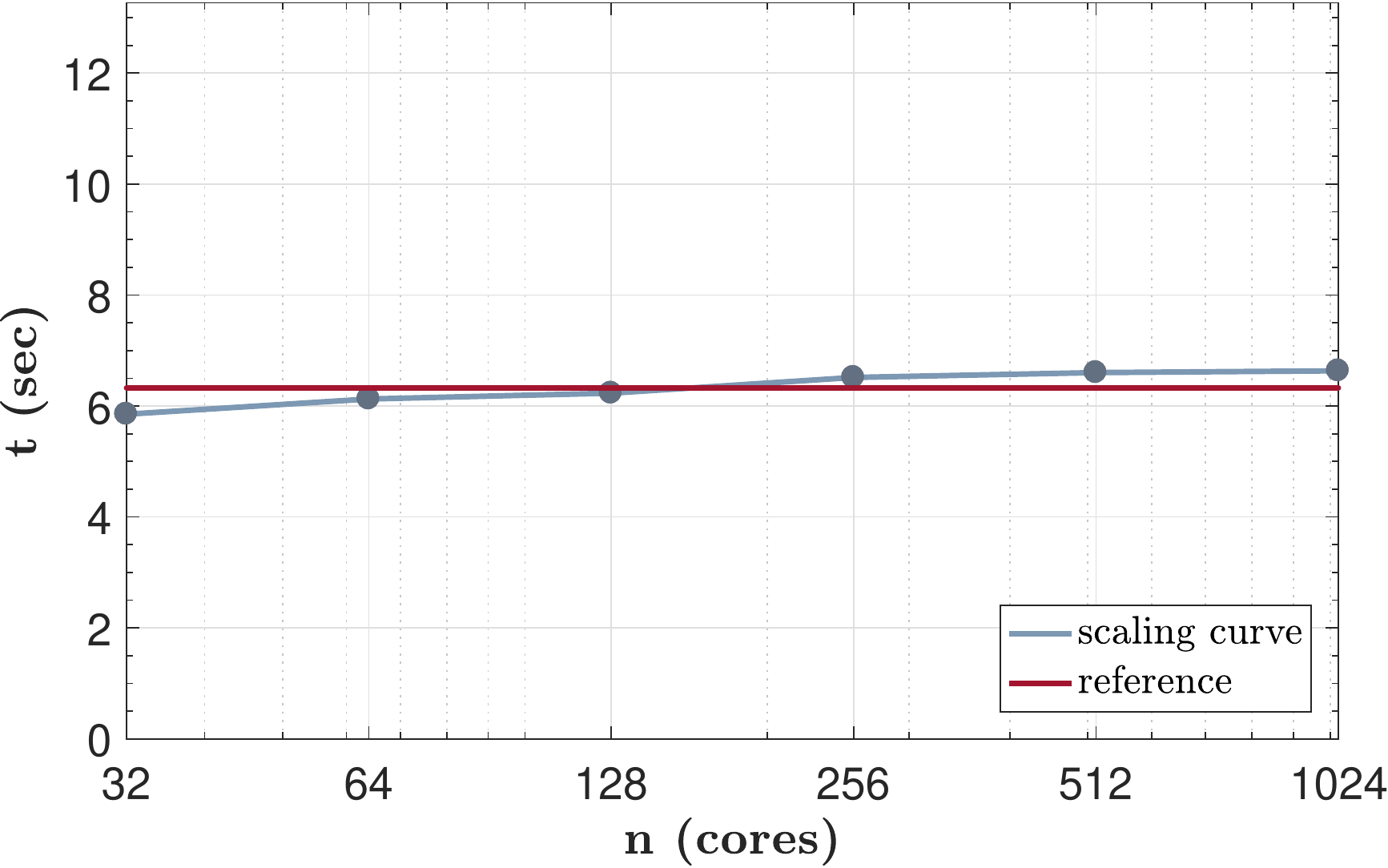}}\qquad
	\subfloat[$2$ million cells per subdomain.]{ \label{fig:ws.ho.3d.b}\includegraphics[width=0.47\textwidth]{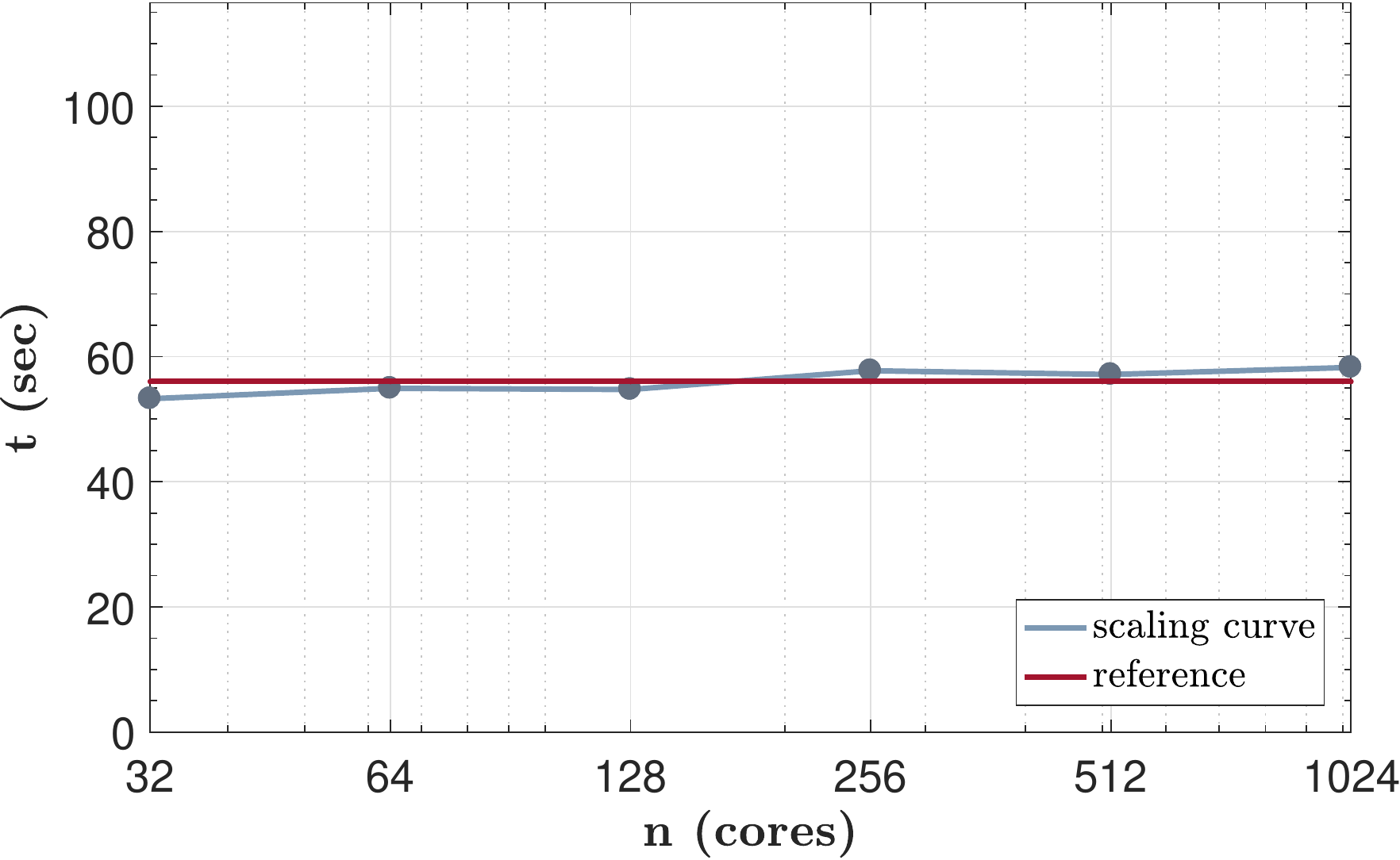}}
	\vspace{15mm}
	\caption{
		%\footnotesize
        Weak scaling curves with homogeneous permeability and 
		$\bar{H} = H$ (see Table \ref{table:ws.ho.3d}). 
	}
	\label{fig:ws.ho.3d}
\end{figure}

\begin{figure}[h!]
	\centering
	\setlength{\abovecaptionskip}{-35pt}
	\subfloat[$262$ thousand cells per subdomain.]{ \label{fig:ws.ho4.3d.a}\includegraphics[width=0.46\textwidth]{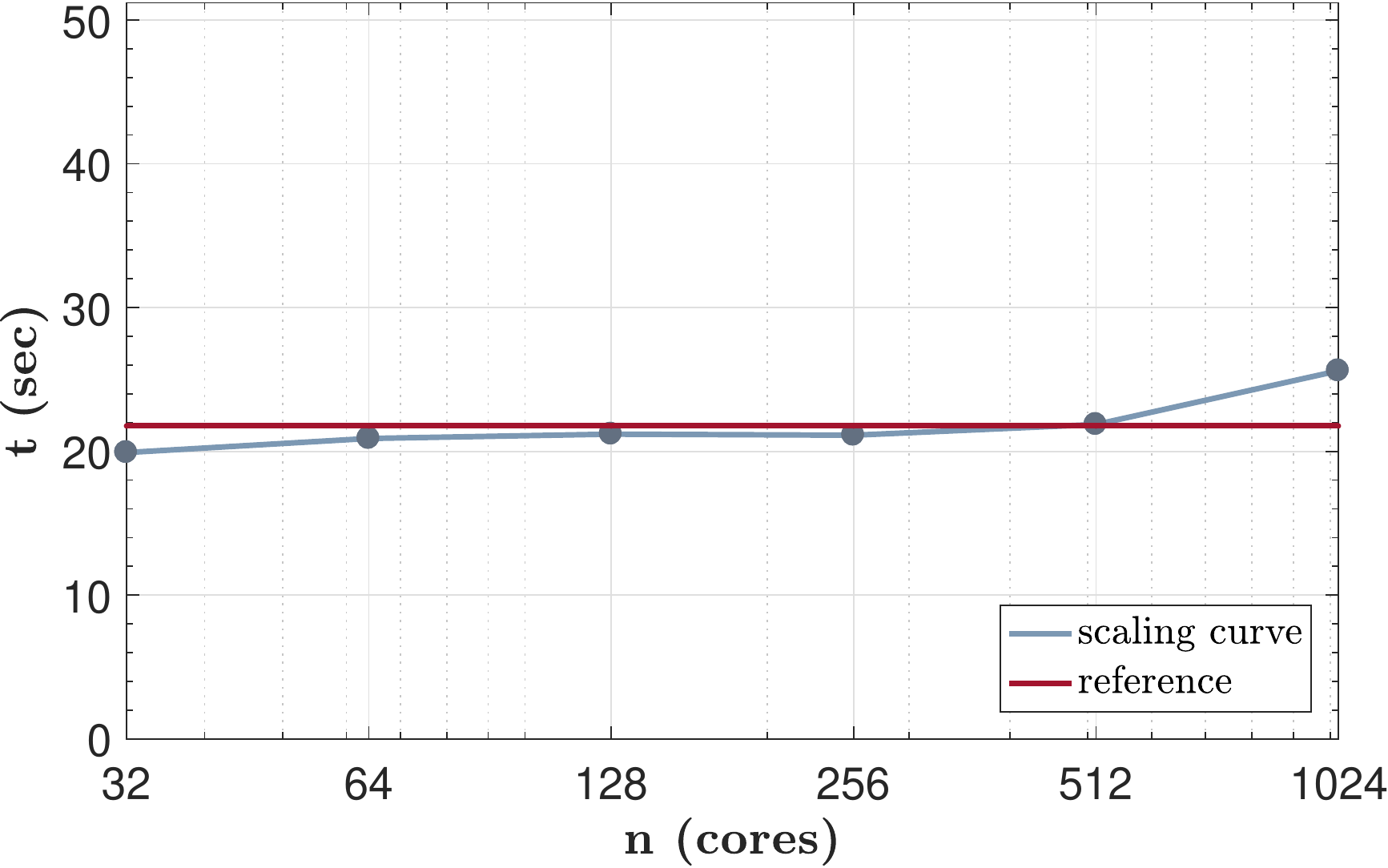}}\qquad
	\subfloat[$2$ million cells per subdomain.]{ \label{fig:ws.ho4.3d.b}\includegraphics[width=0.47\textwidth]{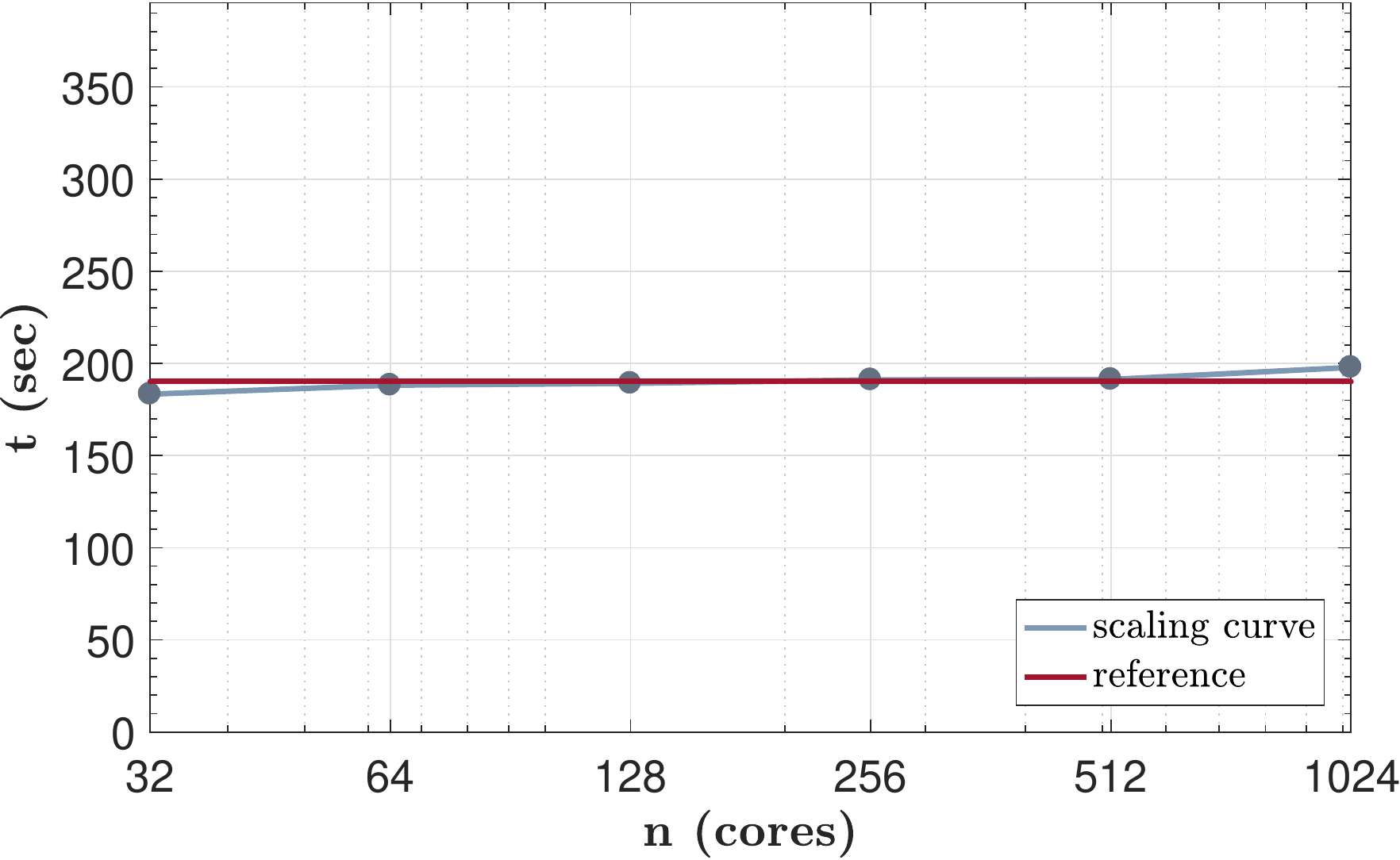}}
	\vspace{15mm}
	\caption{
		%\footnotesize
        Weak scaling curves with homogeneous permeability and 
		$\bar{H} = H/4$ (see Table \ref{table:ws.ho4.3d}). 
	}
	\label{fig:ws.ho4.3d}
\end{figure}

In Figures \ref{fig:ws.ho.3d.a} and \ref{fig:ws.ho.3d.b} we have the scaling
curves with a fixed number of subdomain cells of $ 64 \times 64 \times 64 $
(approximately 262 thousand cells) and $ 128 \times 128 \times 128$ 
(approximately 2 million cells) respectively, computed with $\bar{H} = H$.  Figures
\ref{fig:ws.ho4.3d.a} and \ref{fig:ws.ho4.3d.b}  report the same weak scaling
experiments, but now with $\bar{H} = H/4$.  The figures show that the
computational time of our simulations remains practically constant, which again
shows an outstanding parallel performance, at least up to 1024 processing
cores. 

%\pagebreak

In Tables \ref{table:ws.ho.3d} and \ref{table:ws.ho4.3d}
we have the total time and the individual times of 
the computation of MMBFs, 
the time spent on the interface problem and the message exchange time between cores.
We can see that the total time is essentially constant. 
The time for the exchange of information and the time for the
interface problem are very small, resulting in excellent parallel performance.
There is a slight increase in computational time for the largest cases ($1024$
processing cores), which is related to the computation of the MMBFs by
iterative methods.
%
%There is an increase with the number of cores and 
%on the larger cases, with $ 128 \times 128 \times 128 $ cells.
%However, this increase in the total time is related to the computation of the MMBFs.
%This time could be reduced by using a direct solver for the
%MMBFs.
An important feature of this implementation is the small computational time 
of the interface problem compared to the local problems, and
this is clearly seen in Table \ref{table:ws.ho4.3d}. 
The size and the quantity of interface linear systems to be solved at each level
depends directly on the number of $\bar H$ partitions at the
interface. The problem with $\bar{H} = H$, which translates to one $\bar H$
partition at each subdomain interface, results in an interface linear system of
size $2 \times 2$ in level $\ell = \mathcal{L}$, doubling its size for each previous
level. These systems are small enough to be solved efficiently by a direct
solver based on LU decomposition. 
Also, each set of interface linear systems are solved simultaneously 
within each core, which accelerates the mixed multiscale method.
Now, the total number of interface problems to be solved on all levels 
depends on both the number of $\bar{H}$ partitions and the number of cores.
This means that the more levels we have, more time we are going to spend in 
the computation of the interface problems. This behavior is reported in 
Tables \ref{table:ws.ho.3d} and \ref{table:ws.ho4.3d}, for a fixed $\bar{H}$.
%It specially increases in the problems with an intermediate coarse scale
%$\bar{H} = H/4$, since the size and quantity of interface problems also
%increases. 
%However, they still remain negligible compared to the total time.
% 
Despite this increase, the time spent on interface calculations is still
negligible compare to the total time in the examples considered here.

\begin{table}[h]
\centering
\caption{
Weak scaling times for homogenous problem with $262$ thousand cells
(top table, see Figure \ref{fig:ws.ho.3d.a}) and $2$ million cells (bottom table, see Figure \ref{fig:ws.ho.3d.b}) 
per subdomain. For these problems we considered
$\bar{H} = H$.  
}
\sizefonte
%\begin{tabularx}{\textwidth}{c|c|c|c|c|c}
\begin{tabular}{c|c|c|c|c|c}
\hline \hline
\multicolumn{6}{c}{\multirow{2}{*}{\textbf{Weak Scaling (homogeneous permeability -  $\bar{H} = H$)}}} \\ 
\multicolumn{6}{c}{} \\ 
\hline \hline
\multicolumn{6}{c}{\multirow{2}{*}{\textbf{262 thousand cells per subd.}}} \\ 
\multicolumn{6}{c}{} \\ 
\hline
Cores  &  MMBFs Time (s) & INTRF Time (s) & MPI Time (s)
& Total Time (s) & $\%$ Avg. dev.\\ \hline
32      & 4.21 & 0.0012 & 0.0001 & 5.85 & 6.60 \\
64      & 4.47 & 0.0014 & 0.0002 & 6.13 & 2.22 \\
128     & 4.67 & 0.0058 & 0.0003 & 6.23 & 0.54 \\
256     & 4.83 & 0.0022 & 0.0004 & 6.51 & 3.98 \\
512     & 5.17 & 0.0112 & 0.0007 & 6.60 & 5.38 \\ 
1024    & 5.79 & 0.0922 & 0.0760 & 6.99 & 9.48 \\ \hline
\end{tabular}
%\begin{tabularx}{\textwidth}{c|c|c|c|c|c}
\begin{tabular}{c|c|c|c|c|c}
\multicolumn{6}{c}{\multirow{2}{*}{\textbf{2 million cells per subd. }}}  \\
\multicolumn{6}{c}{} \\
\hline
Cores  &  MMBFs Time (s) & INTRF Time (s) & MPI Time (s)
& Total Time (s) & $\%$ Avg. dev.\\ \hline
32     & 38.87  & 0.0034 & 0.0001 & 53.26 & 4.14 \\
64     & 42.24  & 0.0036 & 0.0002 & 54.93 & 1.14  \\
128    & 43.42  & 0.0039 & 0.0004 & 54.74 & 1.48 \\
256    & 44.21  & 0.0047 & 0.0005 & 57.74 & 3.92 \\
512    & 47.42  & 0.0070 & 0.0008 & 57.14 & 2.84 \\ 
1024   & 53.32  & 0.0893 & 0.0757 & 58.69 & 4.65 \\ \hline
\end{tabular}
\label{table:ws.ho.3d}
\end{table}

\begin{table}[h]
\centering
\caption{
Weak scaling times for homogeneous problem with $262$ thousand cells (top
table, see Figure \ref{fig:ws.ho4.3d.a}) and $2$ million cells (bottom table, see Figure \ref{fig:ws.ho4.3d.b}) 
per subdomain. For these problems we considered
$\bar{H} =H/4$.  
}
\sizefonte
%\begin{tabularx}{\textwidth}{c|c|c|c|c|c}
\begin{tabular}{c|c|c|c|c|c}
\hline \hline
\multicolumn{6}{c}{\multirow{2}{*}{\textbf{Weak Scaling (homogeneous permeability - $\bar{H} = H/4$)}}}\\ 
\multicolumn{6}{c}{}\\ 
\hline \hline
\multicolumn{6}{c}{\multirow{2}{*}{\textbf{262 thousand cells per subd.}}}  \\ 
\multicolumn{6}{c}{}  \\ \hline
Cores &  MMBFs Time (s)  & INTRF Time (s) & MPI
Time (s) & Total Time (s) & $\%$ Avg. dev.\\ \hline
32    & 15.08 & 0.0064 & 0.0006 & 19.91  & 5.16 \\
64    & 16.05 & 0.0120 & 0.0023 & 20.89  & 0.49 \\
128   & 16.15 & 0.0454 & 0.0039 & 21.19  & 0.93 \\
256   & 16.35 & 0.1582 & 0.0110 & 21.12  & 0.59 \\
512   & 17.79 & 0.8692 & 0.0368 & 21.87  & 4.13 \\
1024  & 20.69 & 4.2895 & 0.3349 & 26.05  & 19.29  \\ \hline
\multicolumn{6}{c}{\multirow{2}{*}{\textbf{2 million cells per subd. }} } \\
\multicolumn{6}{c}{} \\ \hline
Cores &  MMBFs Time (s) & INTRF Time (s) & MPI Time (s)
& Total Time (s) & $\%$ Avg. dev. \\ \hline
32    & 139.88  & 0.0132 & 0.0010 & 183.30 & 2.81\\
64    & 149.86  & 0.0194 & 0.0024 & 188.18 & 0.22\\
128   & 151.23  & 0.0525 & 0.0043 & 189.17 & 0.30\\
256   & 153.83  & 0.1618 & 0.0119 & 191.05 & 1.30\\
512   & 165.33  & 0.8364 & 0.0355 & 191.28 & 1.42\\
1024  & 188.91  & 4.1823 & 0.4077 & 198.29 & 4.25 \\ \hline
\end{tabular}
\label{table:ws.ho4.3d}
\end{table}

\clearpage

\subsection{Heterogeneous scalability study}
\label{hetero}

\begin{figure}[h!]
	\centering
	\includegraphics[width=0.49\textwidth]{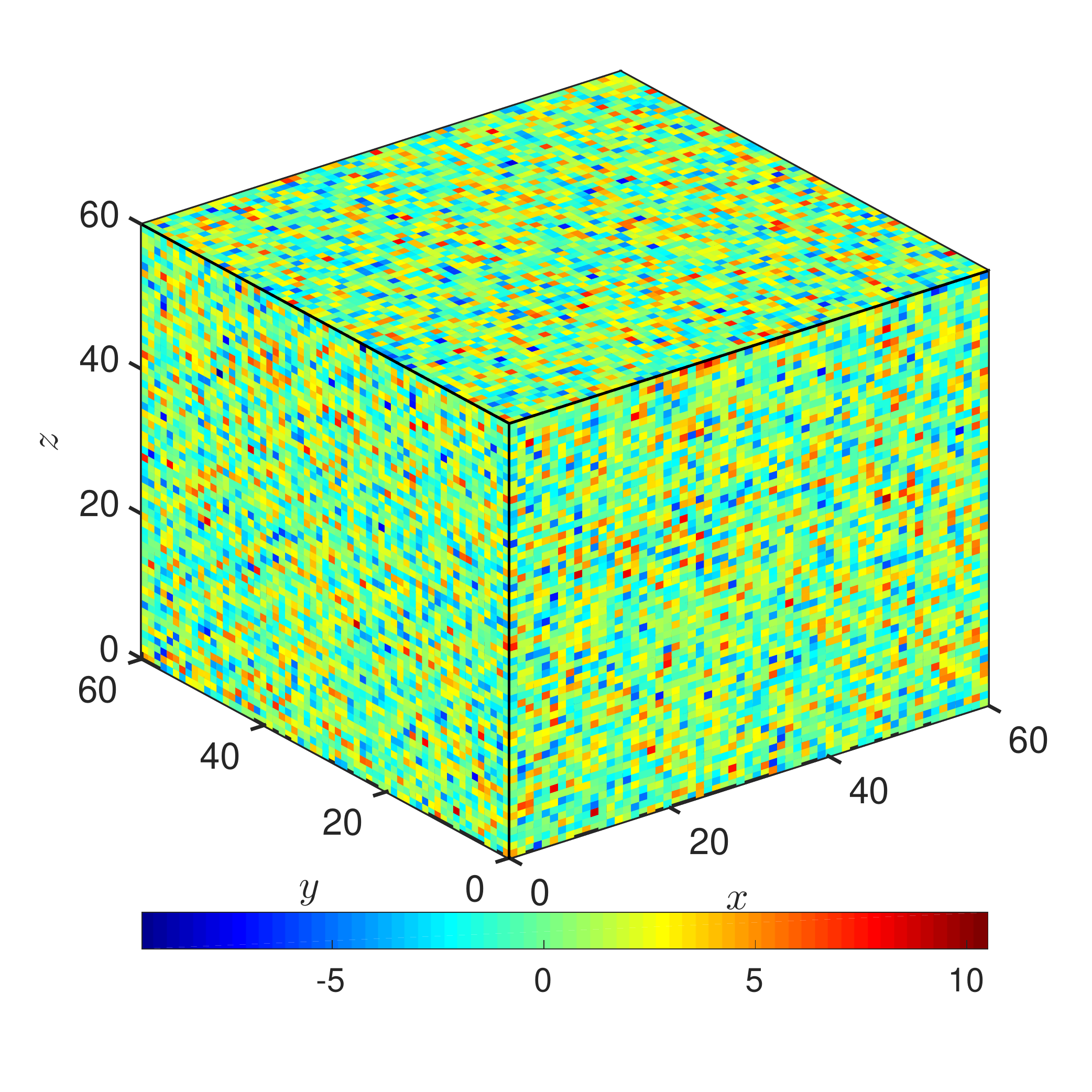}
	\includegraphics[width=0.49\textwidth]{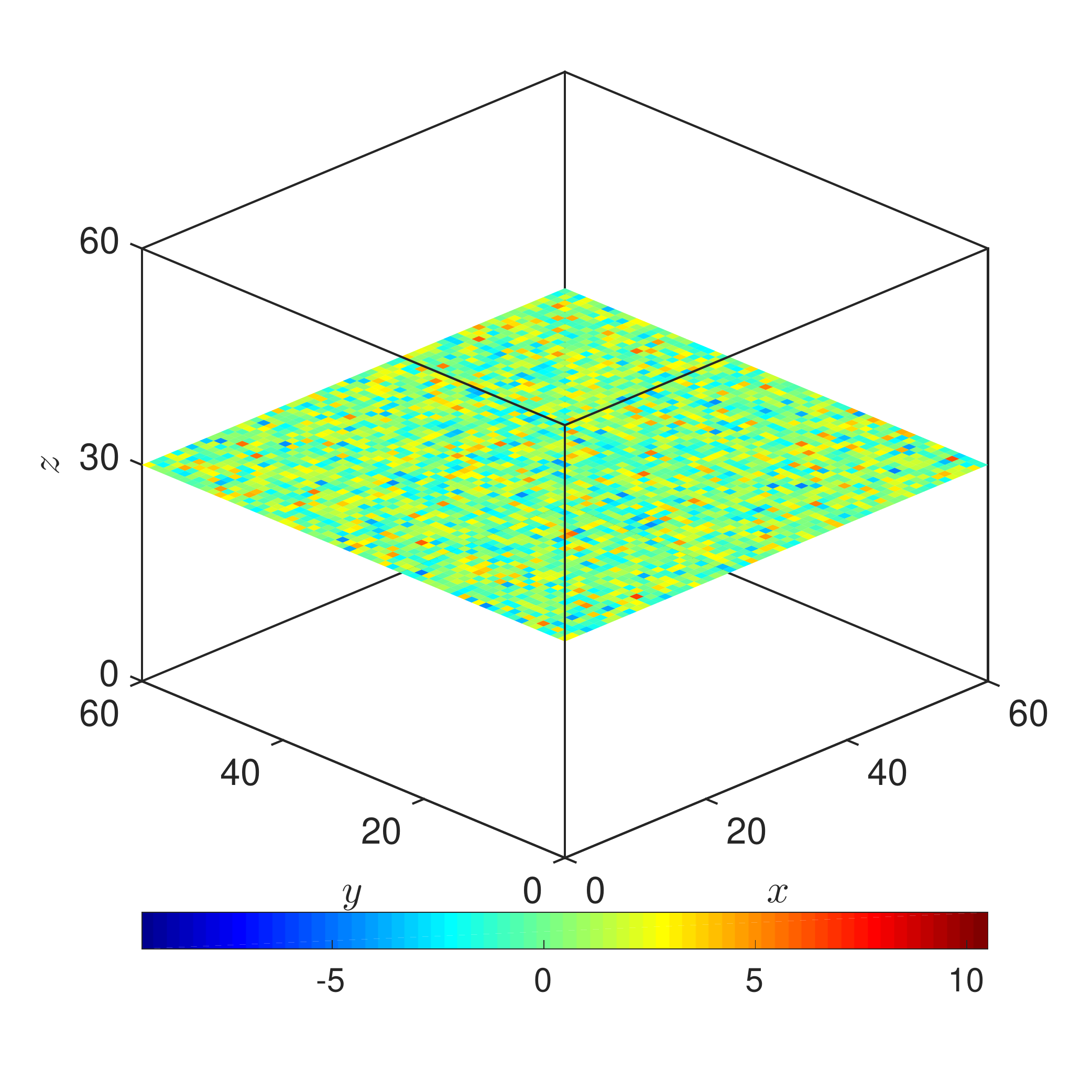}
    \caption{The left picture shows the three-dimensional heterogeneous permeability with a
    contrast of $10^8$ given by \eqref{eq:num.log_normal}. The right picture
    shows a two-dimensional slice on the middle of the $z$-direction of the
    permeability on the left.
    } 
    \label{fig:het_perm}
\end{figure}

For the heterogeneous permeability field, we
use a log-normal model for multiscale rock heterogeneity proposed by Glimm and
Sharp \cite{GS91} where the absolute permeability is given by 
\begin{equation}\label{eq:num.log_normal}
K({\bf x}) = K_0 \exp(\omega_{K}\xi({\bf x})),
\end{equation}
where $\xi(\bf x)$ is an independent Gaussian field with
$K_{0} = 1.6487$ and $\omega_{K} = 3.7$   
in order to generate a permeability field with contrast $K_{max}/K_{min} = 10^8$
on a mesh of $60 \times 60 \times 60$ cells. 
For simulations in finer grid resolutions, the same permeability field is used
by projecting it onto the finer grids.
In Figure \ref{fig:het_perm} we illustrate the three-dimensional heterogeneous 
absolute permeability used for the scaling experiments, and a two-dimensional 
slice in the middle of the $z$-direction that was used for the accuracy experiments. 
%
%\clearpage

\subsubsection{Strong scaling} 

\begin{figure}[h!]
\centering
\setlength{\abovecaptionskip}{-35pt}
\subfloat[$134$ million cells.]{\label{fig:ss.he.3d.a}\includegraphics[width=0.47\textwidth]{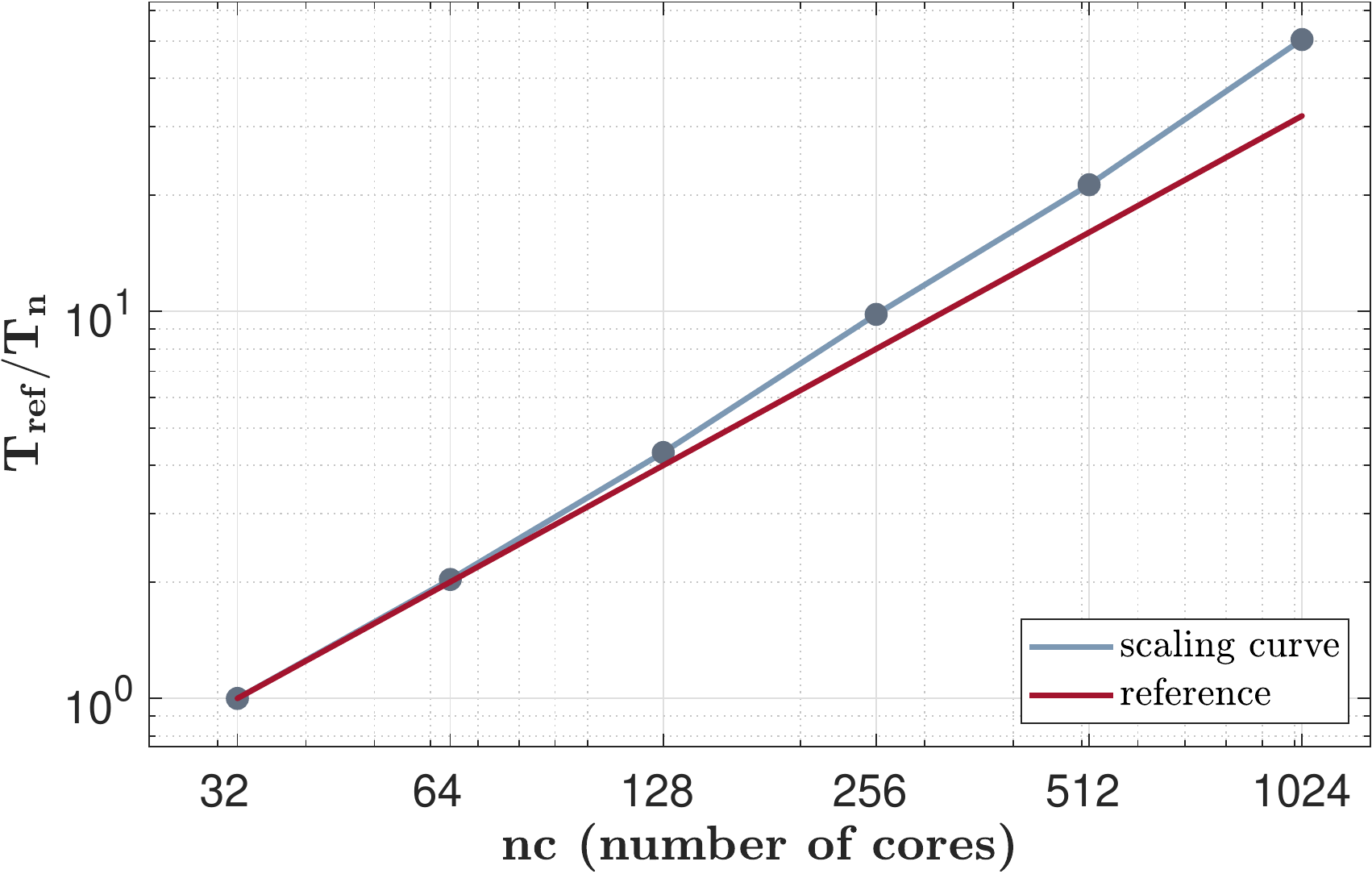}}\qquad
\subfloat[$1$ billion cells.]{\label{fig:ss.he.3d.b}\includegraphics[width=0.46\textwidth]{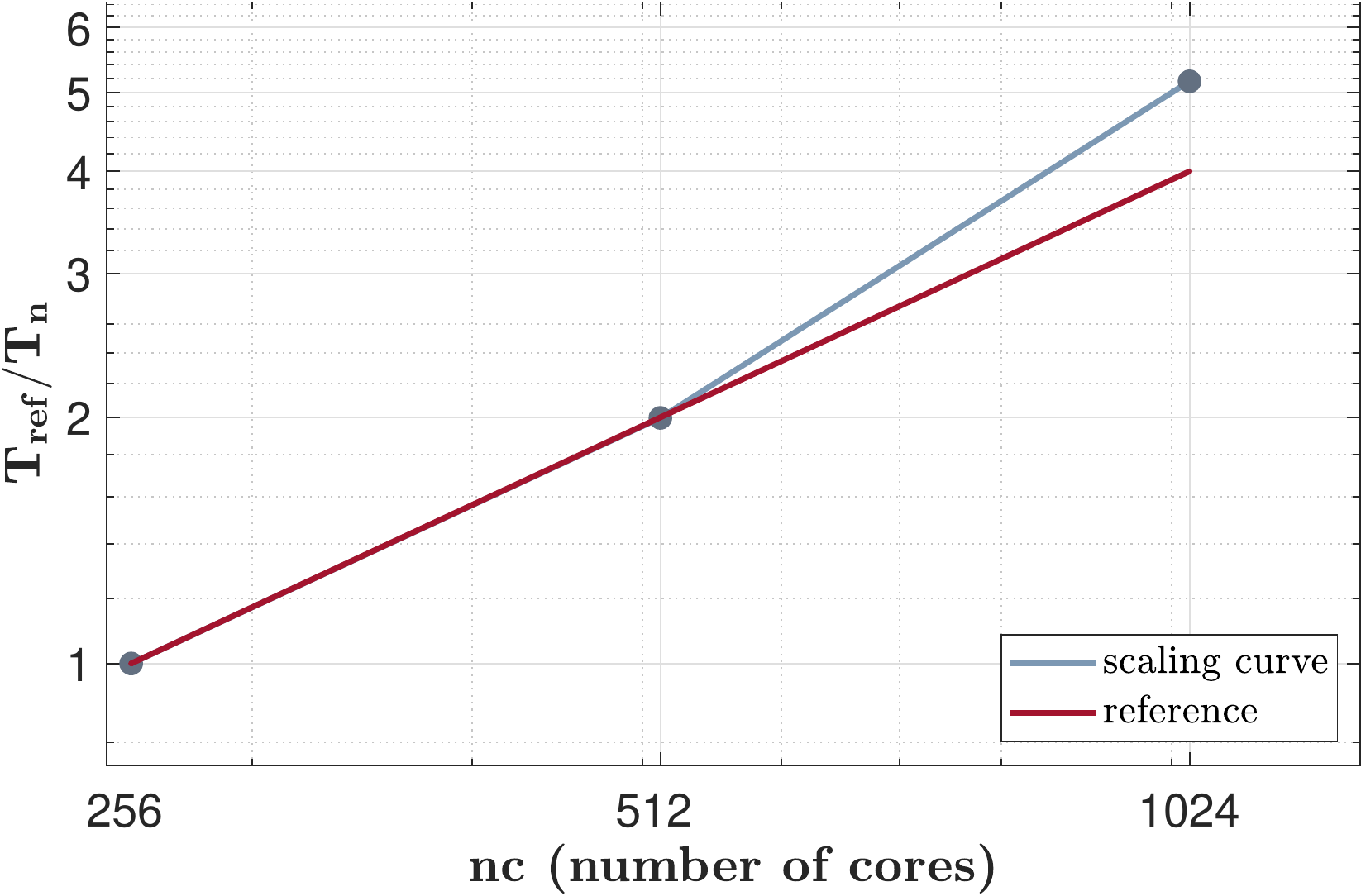}}
\vspace{15mm}
\caption{
Strong scaling curves with heterogeneous log-normal permeability and  
$\bar{H} = H/4$ (see Table \ref{table:ss.he.3d}).
}
\label{fig:ss.he.3d}
\end{figure}

\begin{figure}[h!]
\centering
\setlength{\abovecaptionskip}{-35pt}
\subfloat[$134$ million cells.]{\label{fig:ss.he4.3d.a}\includegraphics[width=0.47\textwidth]{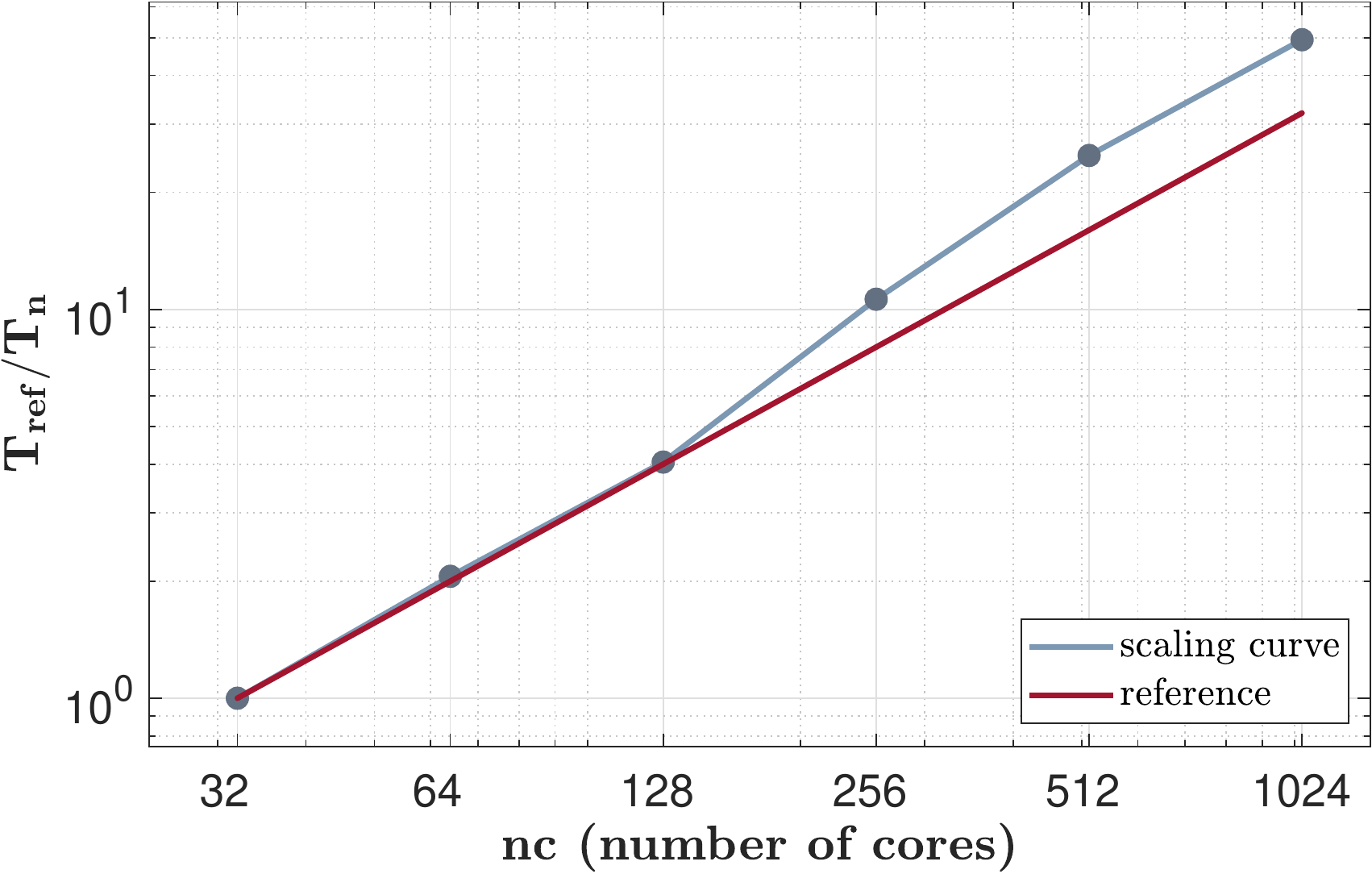}}\qquad
\subfloat[$1$ billion cells.]{\label{fig:ss.he4.3d.b}\includegraphics[width=0.46\textwidth]{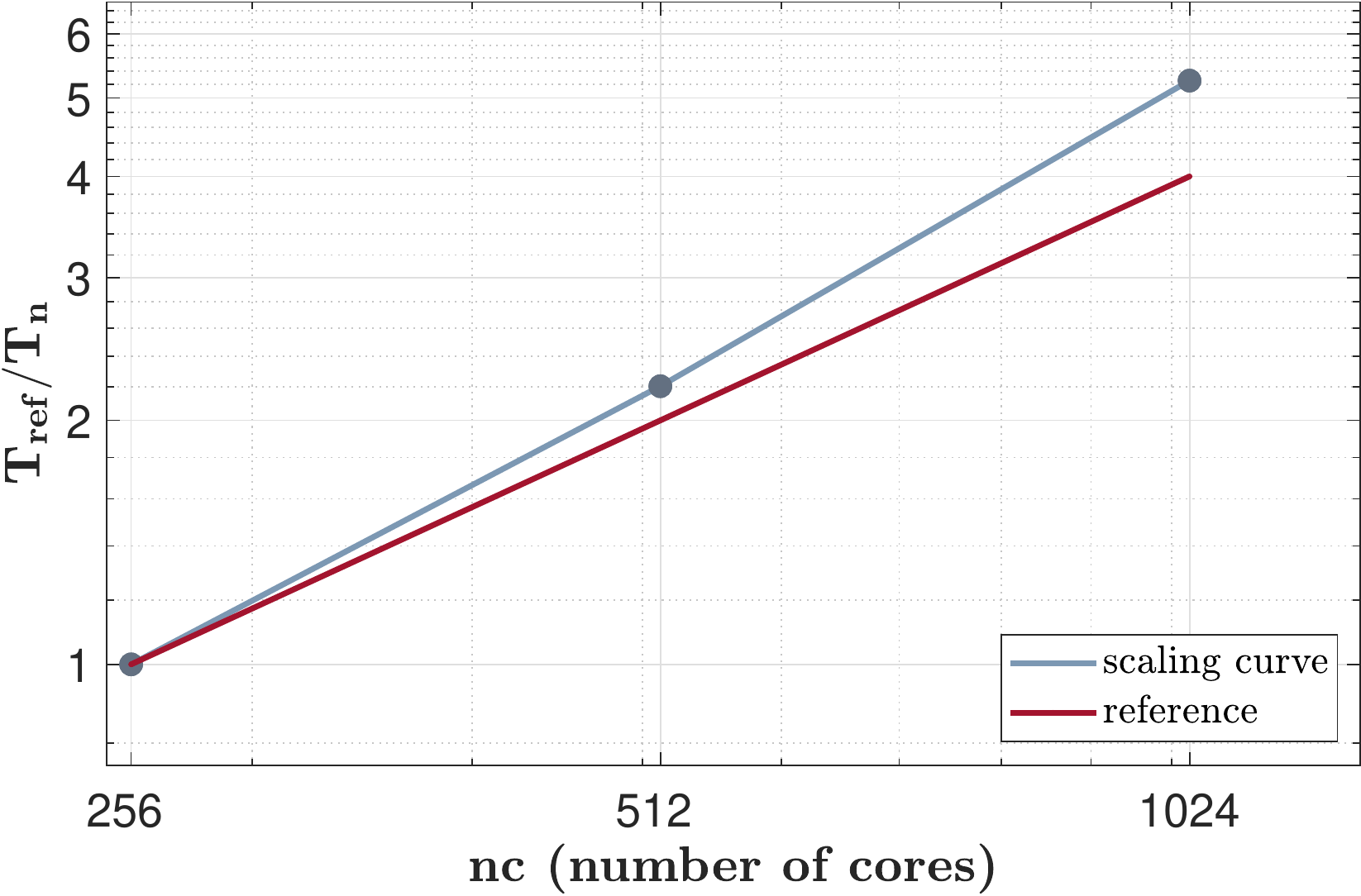}}
\vspace{15mm}
\caption{
%\footnotesize
Strong scaling curves with heterogeneous
log-normal permeability and 
$\bar{H} = H/4$ (see Table \ref{table:ss.he4.3d}).} 
\label{fig:ss.he4.3d}
\end{figure}

For the heterogeneous experiments we considered the same conditions as
homogeneous strong scaling experiments.  In Figures \ref{fig:ss.he.3d.a} and
\ref{fig:ss.he.3d.b}, and Table \ref{table:ss.he.3d} we present the scaling
curves and computational times for the $\bar{H} = H$ problem with total $134$
million and $1$ billion cells, respectively.  In Figures \ref{fig:ss.he4.3d.a}
and \ref{fig:ss.he4.3d.b}, and Table \ref{table:ss.he4.3d} 
we show the experiments considering $\bar{H} = H/4$.
We can see that our simulations again exhibit excellent performance.
The times for the interface problem remain negligible.

\begin{table}[h!]
\centering
\caption{
Strong scaling for heterogeneous problem with $134$ million
cells (top table, see Figure \ref{fig:ss.he.3d.a}) and with $1$ billion
cells (bottom table, see Figure \ref{fig:ss.he.3d.b}). 
For these problems we considered $\bar{H} = H$. 
}
\sizefonte
%\begin{tabularx}{\textwidth}{c|c|Y|Y|Y}
\begin{tabular}{c|c|c|c|c|c}
\hline \hline
\multicolumn{6}{c}{\multirow{2}{*}{\textbf{Strong Scaling (heterogeneous 
$\bar{H} = H$)}}} \\ 
\multicolumn{6}{c}{} \\ 
\hline \hline
\multicolumn{6}{c}{\multirow{2}{*}{\textbf{134 million cells}}} \\ 
\multicolumn{6}{c}{} \\ 
\hline
Cores &  MMBFs Time (s) &  INTRF Time (s) & MPI Time (s) & Total Time (s) & $\%$ decrease (Total) \\ \hline
32    & 303.99  & 0.0225 & 0.0077 & 398.33  &         \\
64    & 143.84  & 0.0055 & 0.0002 & 195.96  & 50.80   \\
128   &  64.88  & 0.0659 & 0.0004 &  91.81  & 53.15   \\
256   &  24.48  & 0.0037 & 0.0005 &  40.21  & 56.21   \\
512   &  11.10  & 0.0060 & 0.0007 &  18.21  & 54.71   \\
1024  &   3.85  & 0.0123 & 0.0023 &   7.48  & 58.91   \\ \hline
\end{tabular}
%\begin{tabularx}{\textwidth}{c|c|Y|Y|Y}
\begin{tabular}{c|c|c|c|c|c}
\multicolumn{6}{c}{\multirow{2}{*}{\textbf{1 billion cells}}} \\ 
\multicolumn{6}{c}{} \\ 
\hline
Cores &  MMBFs Time (s) & INTRF Time (s) & MPI Time (s) & Total Time (s) & $\%$ decrease (Total) \\ \hline
256   & 207.53 & 0.0110 & 0.0007 & 328.40 &        \\
512   &  99.21 & 0.0121 & 0.0008 & 164.05 & 50.05   \\
1024  &  33.18 & 0.0156 & 0.0024 &  63.37 & 61.37   \\ \hline
\end{tabular}
\label{table:ss.he.3d}
\end{table}

% contrast 10+8
\begin{table}[h!]
\centering
\caption{
Strong scaling times for heterogeneous problem with $134$ million
cells (top table, see Figure \ref{fig:ss.he4.3d.a}) and with $1$ billion
cells (bottom table, see Figure \ref{fig:ss.he4.3d.b}). 
For these problems we considered $\bar{H} = H/4$. 
}
\sizefonte
%\begin{tabularx}{\textwidth}{c|c|Y|Y|Y}
\begin{tabular}{c|c|c|c|c|c}
\hline \hline
\multicolumn{6}{c}{\multirow{2}{*}{\textbf{Strong Scaling (heterogeneous 
$\bar{H} = H/4$)}}} \\ 
\multicolumn{6}{c}{} \\ 
\hline \hline
\multicolumn{6}{c}{\multirow{2}{*}{\textbf{134 million cells}}} \\ 
\multicolumn{6}{c}{} \\ 
\hline
Cores &  MMBFs Time (s) & INTRF Time (s) & MPI Time (s) & Total Time (s) & $\%$ decrease (Total) \\ \hline
32    & 1102.10  & 0.1650 & 0.0540 & 1389.00 &      \\
64    &  507.26  & 0.0255 & 0.0024 &  674.43 & 51.44      \\
128   &  233.85  & 0.0558 & 0.0036 &  342.57 & 49.21     \\
256   &   86.87  & 0.1587 & 0.0123 &  130.41 & 61.93     \\
512   &   39.98  & 0.8985 & 0.0350 &   55.24 & 57.64     \\
1024  &   15.80  & 4.2353 & 0.1949 &   27.51 & 50.21     \\ \hline
\end{tabular}
%\begin{tabularx}{\textwidth}{c|c|Y|Y|Y}
\begin{tabular}{c|c|c|c|c|c}
\multicolumn{6}{c}{\multirow{2}{*}{\textbf{1 billion cells}}} \\ 
\multicolumn{6}{c}{} \\ 
\hline
Cores &  MMBFs Time (s) & INTRF Time (s) & MPI Time (s) & Total Time (s) & $\%$ decrease (Total) \\ \hline
256   & 712.82 & 0.1820 & 0.0130 & 1044.30  &     \\
512   & 351.85 & 0.8917 & 0.0337 &  473.37  & 54.67   \\
1024  & 131.19 & 4.0390 & 0.2042 &  198.25  & 58.12   \\ \hline
\end{tabular}
\label{table:ss.he4.3d}
\end{table}

%\newpage
%\clearpage

\subsubsection{Weak scaling} 

For these experiments we considered the same conditions as the previous
experiments of weak scaling but with the same permeability
of the heterogeneous strong scaling study. The heterogeneous 
permeability field is repeated on each subdomain to reproduce the same 
experiment performed for the homogeneous weak scaling. We also use the same 
nondimensionalized boundary conditions.
Notice that the for each set of processors we solve a different 
heterogeneous problem. 
The idea is to keep the same computational effort on each subdomain 
and assess the method behavior especially for the interface problems.

\begin{figure}[h!]
	\centering
	\setlength{\abovecaptionskip}{-35pt}
	\subfloat[$262$ thousand cells per subd.]{ \label{fig:ws.he.3d.a}\includegraphics[width=0.46\textwidth]{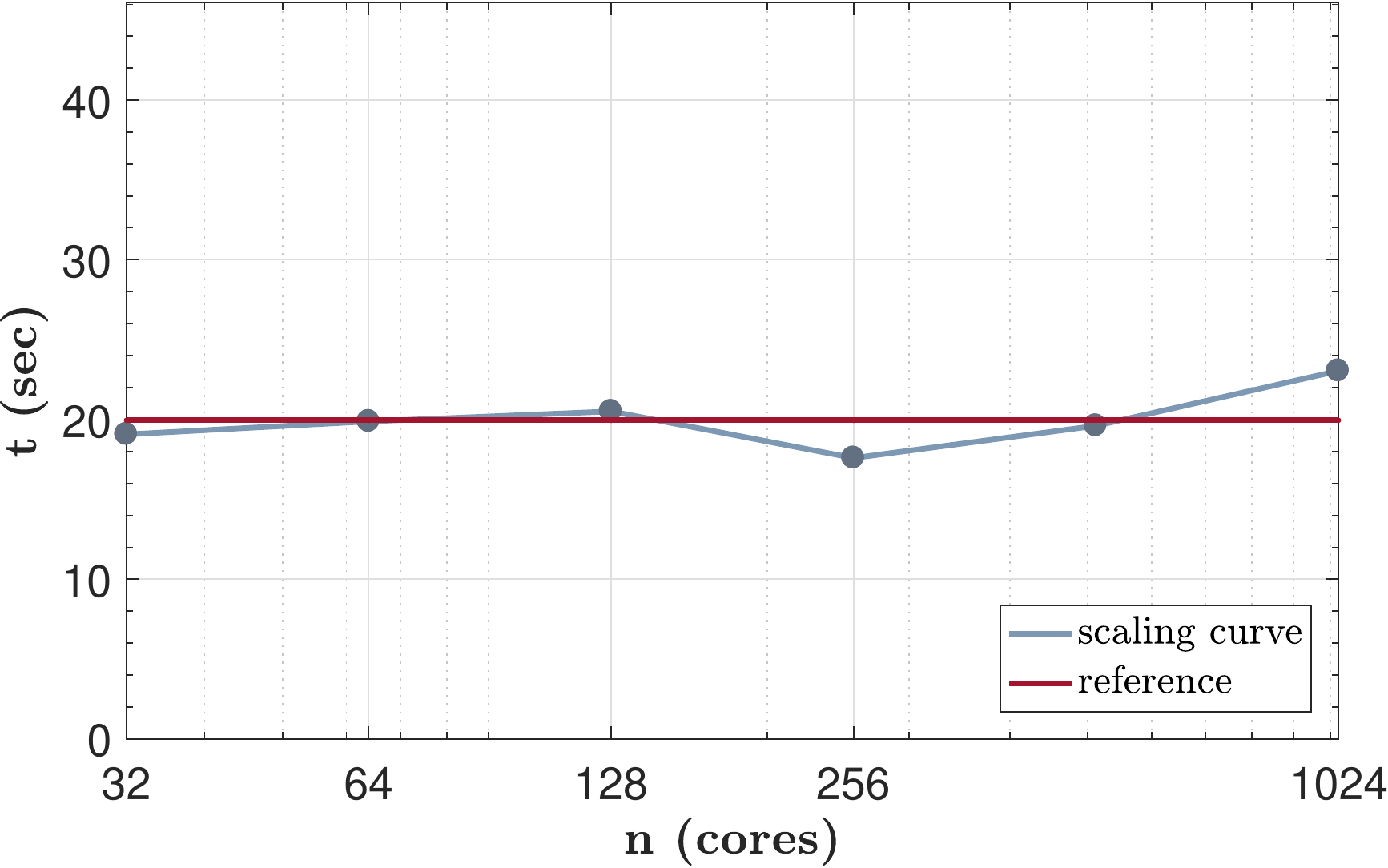}}
	\subfloat[$2$ million cells per subd.]{ \label{fig:ws.he.3d.b}\includegraphics[width=0.47\textwidth]{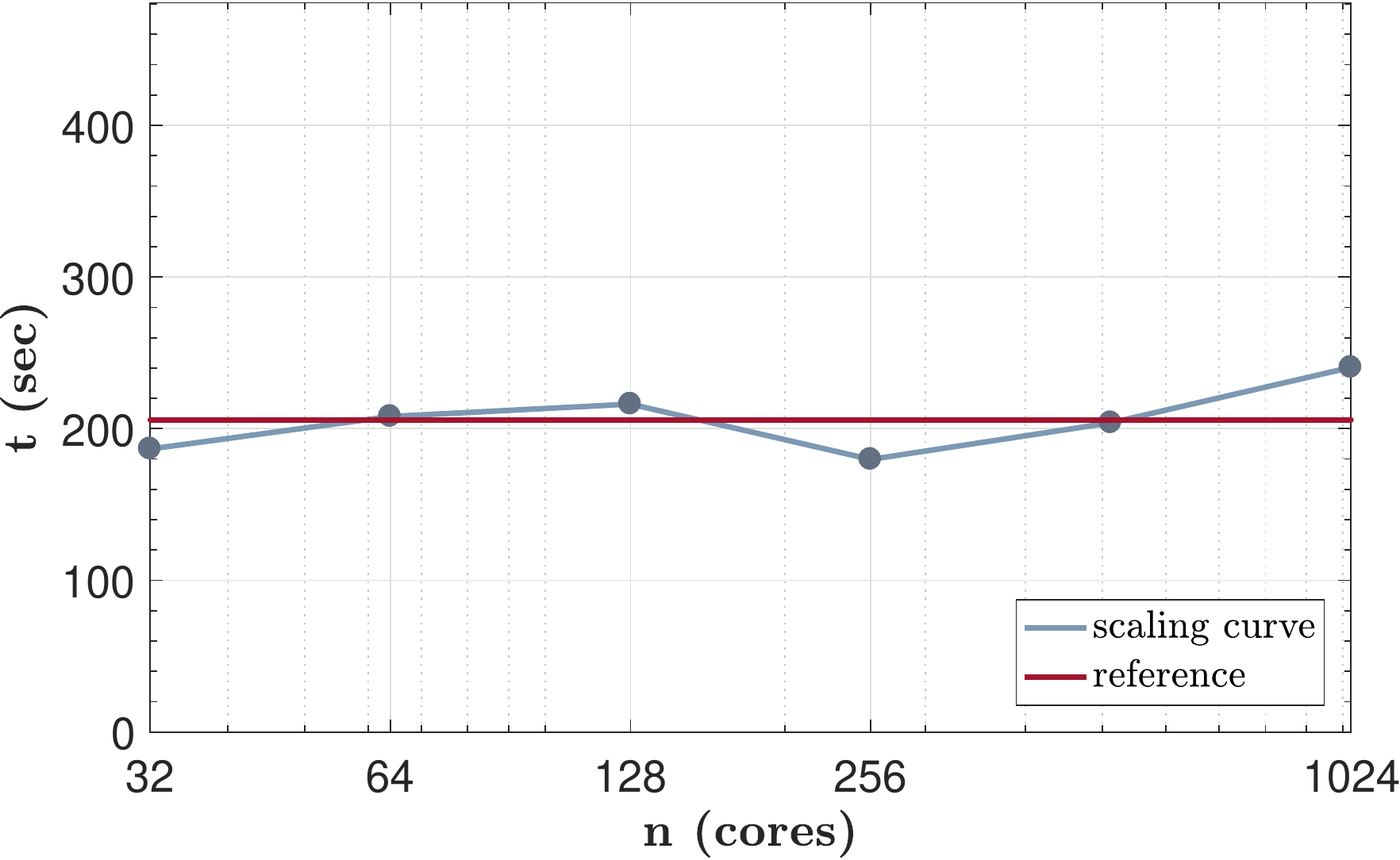}}
	\vspace{15mm}
	\caption{
		%\footnotesize
        Weak scaling curves with heterogeneous log-normal permeability and 
		$\bar{H} = H$ (see Table \ref{table:ws.he.3d}). 
	}
	\label{fig:ws.he.3d}
\end{figure}

\begin{figure}[h!]
	\centering
	\setlength{\abovecaptionskip}{-35pt}
	\subfloat[$262$ thousand cells per subd.]{ \label{fig:ws.he4.3d.a}\includegraphics[width=0.46\textwidth]{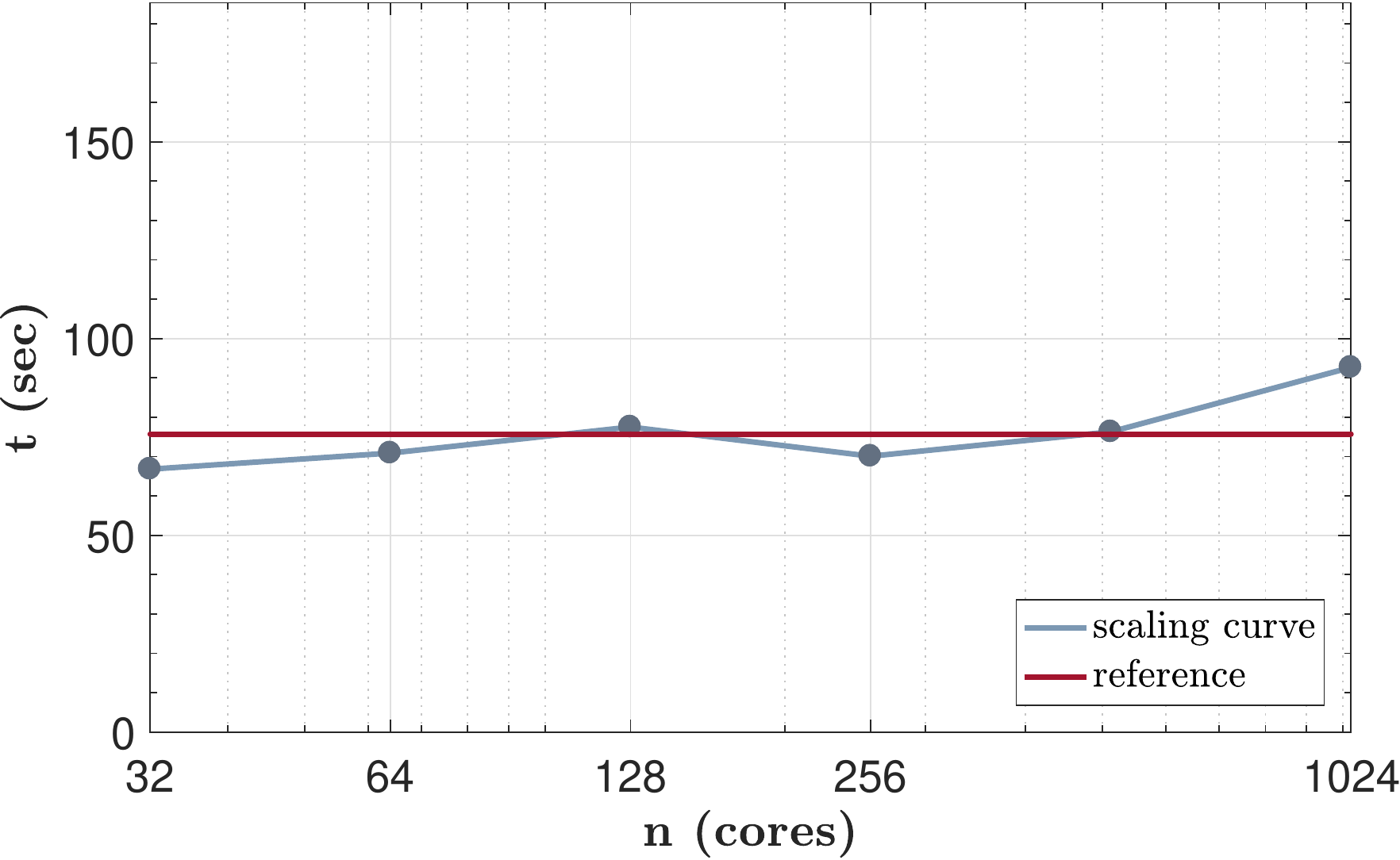}}
	\subfloat[$2$ million cells per subd.]{ \label{fig:ws.he4.3d.b}\includegraphics[width=0.47\textwidth]{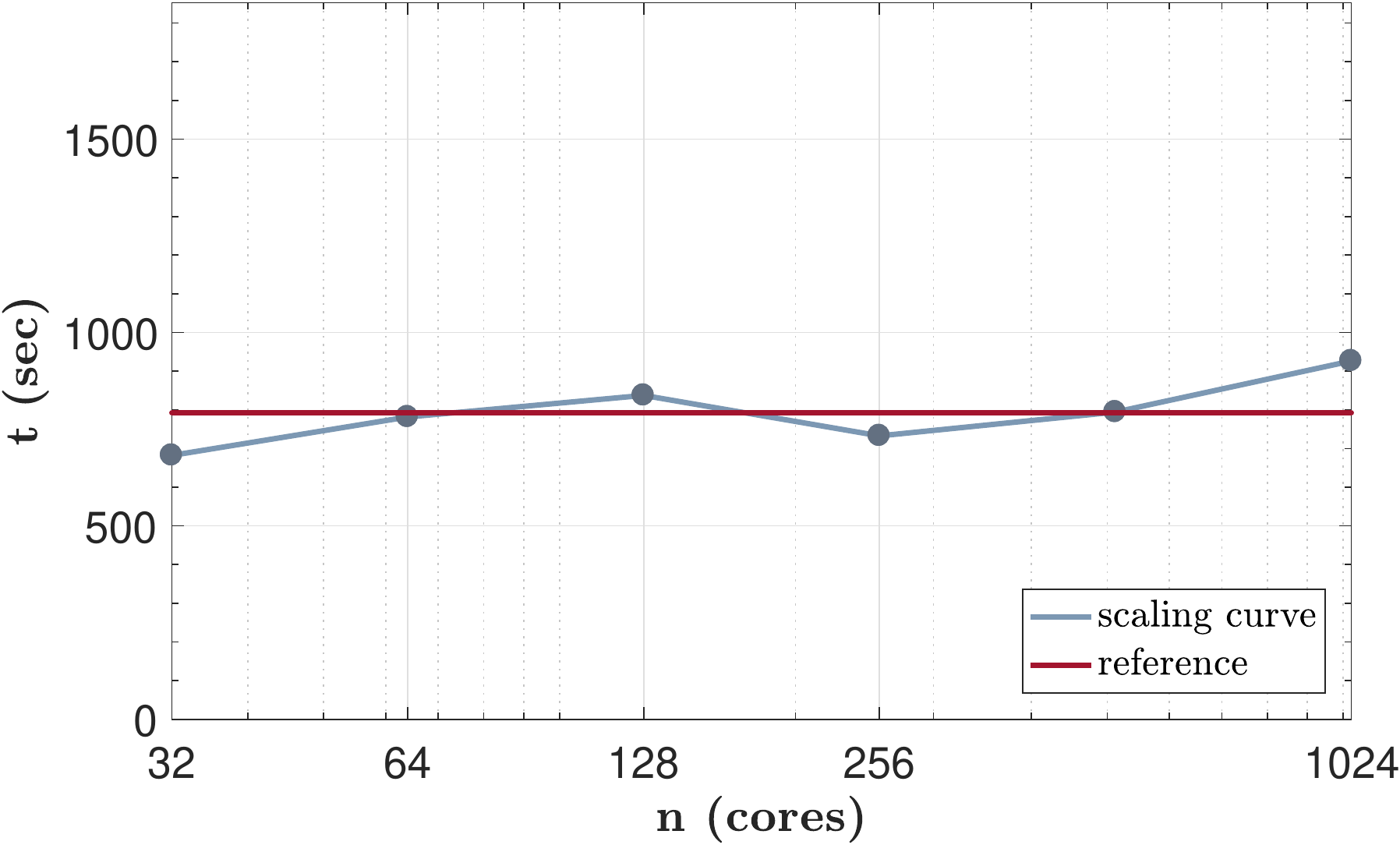}}
	\vspace{15mm}
	\caption{
		%\footnotesize
        Weak scaling curves with heterogeneous log-normal permeability and 
		$\bar{H} = H/4$ (see Table \ref{table:ws.he4.3d}). 
	}
	\label{fig:ws.he4.3d}
\end{figure}

In Figures \ref{fig:ws.he.3d.a} and
\ref{fig:ws.he.3d.b} we have the scaling curves 
for the heterogeneous problem with $\bar{H} = H$ and 
a fixed number of subdomain cells of $262$ thousand and 
$2$ million cells respectively.  
In Figures \ref{fig:ws.he4.3d.a} and \ref{fig:ws.he4.3d.b}  
we show the same weak scaling experiments with $\bar{H} = H/4$.
Tables \ref{table:ws.he.3d} and \ref{table:ws.he4.3d} shows the 
computational times and we can see that the runtime is fairly close. 
The variation in time that we see is related to the fact that we
solve different problems at each point since we repeat the permeability in each
subdomain. Still, the influence of the MMBFs computation 
dominates the total time. 
The computation time of the interface problems continue to have 
little influence on the total time and follows the same increase 
pattern of the homogeneous problem as we consider the same number 
of $\bar H$ partitions and levels.  This shows that the interface 
time is insensible to changes in the permeability.
\begin{table}[h!]
\centering
\caption{
Weak scaling times for heterogeneous problem with $262$ thousand cells
(top table, see Figure \ref{fig:ws.he.3d.a}) and $2$ million cells (bottom table, see Figure \ref{fig:ws.he.3d.b})
per subdomain. For these problems we considered $\bar{H} = H$.
}
\sizefonte
%\begin{tabularx}{\textwidth}{c|c|c|c|c|c}
\begin{tabular}{c|c|c|c|c|c}
\hline \hline
\multicolumn{6}{c}{\multirow{2}{*}{\textbf{Weak Scaling (heterogeneous $\bar{H} = H$)}}} \\ 
\multicolumn{6}{c}{} \\ 
\hline \hline
\multicolumn{6}{c}{\multirow{2}{*}{\textbf{262 thousand cells per subd.}}} \\ 
\multicolumn{6}{c}{} \\ 
\hline
Cores  &  MMBFs Time (s) & INTRF Time (s) & MPI Time (s) 
& Total Time (s) & $\%$ Avg. dev.\\ \hline
32     & 18.22 & 0.0012 & 0.0001 & 19.42 & 4.78   \\
64     & 18.84 & 0.0017 & 0.0002 & 20.36 & 0.17   \\
128    & 19.47 & 0.0021 & 0.0004 & 21.00 & 2.97   \\
256    & 16.67 & 0.0028 & 0.0006 & 17.98 & 11.84  \\
512    & 18.42 & 0.0158 & 0.0051 & 20.16 & 1.15   \\
1024   & 21.93 & 0.0813 & 0.0334 & 23.45 & 14.98  \\ \hline
\end{tabular}
%\begin{tabularx}{\textwidth}{c|c|c|c|c|c}
\begin{tabular}{c|c|c|c|c|c}
\multicolumn{6}{c}{\multirow{2}{*}{\textbf{2 million cells per subd. }}}  \\ 
\multicolumn{6}{c}{} \\ 
\hline
Cores  &  MMBFs Time (s) & INTRF Time (s) & MPI Time (s) 
& Total Time (s) & $\%$ Avg. dev.\\ \hline
32     & 181.89 & 0.0037 & 0.0001 & 186.98 & 9.35  \\
64     & 203.14 & 0.0043 & 0.0002 & 208.53 & 1.10  \\
128    & 211.49 & 0.0046 & 0.0006 & 216.69 & 5.06  \\
256    & 174.90 & 0.0123 & 0.0044 & 180.12 & 12.67 \\
512    & 198.65 & 0.0138 & 0.0015 & 204.43 & 0.89  \\
1024   & 235.34 & 0.0232 & 0.0085 & 240.82 & 16.75 \\ \hline
\end{tabular}
\label{table:ws.he.3d}
\end{table}

\begin{table}[h!]
\centering
\caption{
Weak scaling times for heterogeneous problem with $262$ thousand cells (top
table, see Figure \ref{fig:ws.he4.3d.a}) and $2$ million cells (bottom table, see Figure \ref{fig:ws.he4.3d.b}) 
per subdomain. For these problems we considered $\bar{H} = H/4$. 
}
\sizefonte
%\begin{tabularx}{\textwidth}{c|c|c|c|c|c}
\begin{tabular}{c|c|c|c|c|c}
\hline \hline
\multicolumn{6}{c}{\multirow{2}{*}{\textbf{Weak Scaling (heterogeneous $\bar{H} = H/4$) }}}  \\ 
\multicolumn{6}{c}{} \\ 
\hline \hline
\multicolumn{6}{c}{\multirow{2}{*}{\textbf{262 thousand cells per subd.}}}  \\ 
\multicolumn{6}{c}{} \\ 
\hline
Cores  &  MMBFs Time (s) & INTRF Time (s) & MPI time (s)
& Total Time (s) & $\%$ Avg. dev. \\ \hline
32    & 65.95  & 0.0056 & 0.0007 & 67.18  & 11.66  \\
64    & 69.83  & 0.0132 & 0.0026 & 71.38  & 6.14   \\
128   & 76.43  & 0.0491 & 0.0054 & 78.01  & 2.58   \\
256   & 68.94  & 0.2126 & 0.0545 & 70.47  & 7.33   \\
512   & 73.67  & 1.4458 & 0.0692 & 76.89  & 1.11   \\
1024  & 87.35  & 4.4946 & 0.5147 & 93.20  & 22.55  \\ \hline
\end{tabular}
%\begin{tabularx}{\textwidth}{c|c|c|c|c|c}
\begin{tabular}{c|c|c|c|c|c}
\multicolumn{6}{c}{\multirow{2}{*}{\textbf{2 million cells per subd.}}}  \\ 
\multicolumn{6}{c}{} \\ 
\hline
Cores  &  MMBFs Time (s) & INTRF Time (s) & MPI Time (s)
& Total Time (s) & $\%$ Avg. dev. \\ \hline
32    & 677.54  & 0.0131 & 0.0013 & 682.59 & 13.87  \\
64    & 776.62  & 0.0222 & 0.0033 & 782.04 & 1.32   \\
128   & 832.58  & 0.0552 & 0.0052 & 837.85 & 5.72   \\
256   & 727.69  & 0.2032 & 0.0454 & 733.13 & 7.49   \\
512   & 787.45  & 1.4758 & 0.0782 & 794.79 & 0.29   \\
1024  & 917.05  & 4.2828 & 0.5280 & 926.82 & 16.95  \\ \hline 
\end{tabular}
\label{table:ws.he4.3d}
\end{table}

%

%{\color{blue}
%To the best of our knowledge of current events, there are very few promising
%Weak/Strong scaling results for parallelizable multiscale methods in the
%literature (see, \cite{manea16,manea17,manea19,MT19,ZYSW17,PED18}). 
%For concreteness, in these papers the authors
%consider a small number of processing cores and numerical simulations with at
%most hundreds of million cells, which are comparable to our
%smallest simulations. In addition, in these papers the most extensive
%scalability study for multiscale methods is made for an algebraic extension of
%multiscale methods \cite{manea16,manea17,manea19} that uses the multiscale
%method as a preconditioner to damp low/high frequency modes. It is noteworthy that in
%this current work, we consider very large problems that can have billions of
%discretization cells, motivated by the numerical simulation of subsurface
%flows. On the other hand, one might find the numerical experiments reported in
%those papers are only restrict to 128 million cells for both GPU
%\cite{manea19} and CPU shared-memory architecture up to 20 processing cores on a single
%node \cite{manea16,manea17}. Moreover, it is important to note that although a strong
%scalability up to 13 speedup for 20 processing cores is found in papers \cite{manea16,manea17,manea19}, 
%the results were obtained for problems much smaller than we dealing in this paper.
%
%}

\clearpage

\subsection{Velocity accuracy in scaling studies}\label{sec:3.5.1}

In this section we present a study on the velocity field error to assess the behavior
of the solution as we increase the number of cores. 
In all studies we compute the accuracy of the flux compared to a reference
solution. The problems presented are two-dimensional slices of the
strong and weak scalings three-dimensional experiments reported above.
We compute 
$||\bu -\bu_h||_{L^2(\Omega)}/||\bu||_{L^2(\Omega)}$, 
the $L^2(\Omega)$ relative velocity error norm, 
where $\bu$ is the reference solution obtained by a hybrid mixed finite element  
discretization \cite{raviart1977mixed,RT91} 
with the same AMG solver \cite{AMGtoolbox} used to solve the
resulting linear system for the pressure. 

In Tables \ref{tab:ss_er} and \ref{tab:sshet_er} we present the relative velocity error norms of the
two-dimensional slice of the strong scaling problems with $134$ million cells 
for the homogeneous permeability 
and heterogeneous permeability (see Figure \ref{fig:het_perm}), 
respectively. 
The errors shown are obtained both for $\bar{H} = H$
and $\bar{H} = H/4$. The tables indicate that as the core number increases (and
consequently the number of levels increases) the
approximated solution does not deteriorate. 
\begin{table}[h!]
\centering
%\scriptsize
\caption{
A two-dimensional flux accuracy study for homogeneous problem
with $262$ thousand cells. The problem represents a two-dimensional
slice of the homogeneous strong scaling problem \eqref{table:ss.ho.3d}
and \eqref{table:ss.ho4.3d} with $134$ million cells.
}
\sizefonte
%\begin{tabularx}{\textwidth}{c|c|c|c}
\begin{tabular}{c|c|c|c}
%\hline \hline
%\multicolumn{4}{c}{\multirow{2}{*}{\textbf{Accuracy Study 1.1 (homogeneous permeability)}}} \\
%\multicolumn{4}{c}{} \\
\hline \hline
\multicolumn{4}{c}{\multirow{2}{*}{\textbf{262 thousand cells (homogeneous permeability)}}} \\
\multicolumn{4}{c}{} \\
\hline
\multicolumn{1}{c|}{} & Subd. cells & $\bar{H} = H$ & $\bar{H} = H/4$ \\
Cores  &  ($n_x \times n_y \times n_z$) & $||\bu -\bu_h||_2/||\bu||_2$ & $||\bu -\bu_h||_2/||\bu||_2$ \\ \hline
32   & 128 $\times$ 64 $\times$ 1 & 5.21e-05 & 6.02e-06  \\
64   & 64 $\times$ 64 $\times$  1 & 1.01e-05 & 1.11e-06  \\
128  & 64 $\times$ 32 $\times$  1 & 2.69e-05 & 3.11e-06  \\
256  & 32 $\times$ 32 $\times$  1 & 5.16e-06 & 5.95e-07  \\
512  & 32 $\times$ 16 $\times$  1 & 1.35e-05 & 1.58e-06  \\
1024 & 16 $\times$ 16 $\times$  1 & 2.56e-06 & 3.10e-07   \\ \hline
\end{tabular}
\label{tab:ss_er}
\end{table}

\begin{table}[h!]
\centering
%\scriptsize
\caption{
A two-dimensional flux accuracy study for heterogeneous problem
with $262$ thousand cells. The problem represents a two-dimensional
slice of the heterogeneous strong scaling problem \eqref{table:ss.he.3d}
and \eqref{table:ss.he4.3d} with $134$ million cells.
}
\sizefonte
%\begin{tabularx}{\textwidth}{c|c|c|c}
\begin{tabular}{c|c|c|c}
%\hline \hline
%\multicolumn{5}{c}{\multirow{2}{*}{\textbf{Accuracy Study (heterogeneous permeability)}}} \\
%\multicolumn{5}{c}{} \\
\hline \hline
\multicolumn{4}{c}{\multirow{2}{*}{\textbf{262 thousand cells (heterogeneous permeability)}}} \\
\multicolumn{4}{c}{} \\
\hline
\multicolumn{1}{c|}{} & Subd. cells& $\bar{H} = H$ & $\bar{H} = H/4$ \\
Cores &  ($n_x \times n_y \times n_z$) & $||\bu -\bu_h||_2/||\bu||_2$ & $||\bu -\bu_h||_2/||\bu||_2$ \\ \hline
32   & 128 $\times$ 64 $\times$ 1 & 5.60e-01 & 4.02e-01   \\
64   & 64 $\times$ 64 $\times$  1 & 6.70e-01 & 4.26e-01  \\
128  & 64 $\times$ 32 $\times$  1 & 6.55e-01 & 3.95e-01   \\
256  & 32 $\times$ 32 $\times$  1 & 7.10e-01 & 3.87e-01   \\
512  & 32 $\times$ 16 $\times$  1 & 7.04e-01 & 3.27e-01   \\
1024 & 16 $\times$ 16 $\times$  1 & 7.20e-01 & 2.87e-01    \\ \hline
\end{tabular}
\label{tab:sshet_er}
\end{table}

%\begin{table}[h!]
%\centering
%%\scriptsize
%\caption{
%A two-dimensional flux accuracy study for heterogeneous problem
%with $262$ thousand cells. The problem represents a two-dimensional
%slice of the heterogeneous strong scaling problem \eqref{table:ss.he.3d}
%and \eqref{table:ss.he4.3d} with $134$ million cells.
%}
%\sizefonte
%%\begin{tabularx}{\textwidth}{c|c|c|c}
%\begin{tabular}{c|c|c|c|c}
%%\hline \hline
%%\multicolumn{5}{c}{\multirow{2}{*}{\textbf{Accuracy Study (heterogeneous permeability)}}} \\
%%\multicolumn{5}{c}{} \\
%\hline \hline
%\multicolumn{5}{c}{\multirow{2}{*}{\textbf{262 thousand cells (heterogeneous permeability)}}} \\
%\multicolumn{5}{c}{} \\
%\hline
%\multicolumn{1}{c|}{} & Subd. cells& $\bar{H} = H$ & $\bar{H} = H/4$ \\
%Cores &  ($n_x \times n_y \times n_z$) & $||\bu -\bu_h||_2/||\bu||_2$ & $||\bu -\bu_h||_2/||\bu||_2$ & Error Reduction \\ \hline
%32   & 128 $\times$ 64 $\times$ 1 & 5.60e-01 & 4.02e-01 & 15.8\%  \\
%64   & 64 $\times$ 64 $\times$  1 & 6.70e-01 & 4.26e-01 & 24.4\% \\
%128  & 64 $\times$ 32 $\times$  1 & 6.55e-01 & 3.95e-01 & 26.0\%  \\
%256  & 32 $\times$ 32 $\times$  1 & 7.10e-01 & 3.87e-01 & 32.3\%  \\
%512  & 32 $\times$ 16 $\times$  1 & 7.04e-01 & 3.27e-01 & 37.7\%  \\
%1024 & 16 $\times$ 16 $\times$  1 & 7.20e-01 & 2.87e-01 & 43.3\%   \\ \hline
%\end{tabular}
%\label{tab:sshet_er}
%\end{table}

In Tables \ref{tab:wk_er} and \ref{tab:wkhet_er} we show the relative flux errors 
for a slice of the weak scaling problem with $262$ thousand cells per
subdomain. In Table \ref{tab:wk_er} we present the errors
for the homogeneous permeability problem and 
in Table \ref{tab:wkhet_er} we present the errors for the  
heterogeneous permeability, both for the two-dimensional slice with $64 \times 64
\times 1$ cells per subdomain. 
For the heterogeneous problems the slice shown in Figure \ref{fig:het_perm} is
repeated in each subdomain, using the same strategy as the weak scaling
experiments above.
The errors are again obtained both for $\bar{H} = H$ and $\bar{H} = H/4$. 
The tables indicate that as the core number increases 
there are no loss of accuracy in the approximated solution.  
\begin{table}[h!]
\centering
%\scriptsize
\caption{
A two-dimensional flux accuracy study $1.2$ for homogeneous problem
with $4096$ cells. The problem represents a two-dimensional
slice of the homogeneous weak scaling problem \eqref{table:ws.ho.3d}
and \eqref{table:ws.ho4.3d} with $262$ thousand cells per subdomain.
}
\sizefonte
%\begin{tabularx}{\textwidth}{c|c|c|c|c|c}
\begin{tabular}{c|c|c|c}
%\hline \hline
%\multicolumn{4}{c}{\multirow{2}{*}{\textbf{Accuracy Study 1.2 (homogeneous permeability)}}} \\
%\multicolumn{4}{c}{} \\
\hline \hline
\multicolumn{4}{c}{\multirow{2}{*}{\textbf{4096 cells per subd. (homogeneous permeability)}}} \\
\multicolumn{4}{c}{} \\
\hline
\multicolumn{1}{c|}{} & Subd. cells & $\bar{H} = H$ & $\bar{H} = H/4$ \\
Cores  & ($n_x \times n_y \times n_z$) & $||\bu -\bu_h||_2/||\bu||_2$ & $||\bu -\bu_h||_2/||\bu||_2$ \\ \hline
32   & 64 $\times$ 64 $\times$ 1 & 9.87e-06 & 1.11e-06  \\
64   & 64 $\times$ 64 $\times$ 1 & 1.01e-05 & 1.11e-06  \\
128  & 64 $\times$ 64 $\times$ 1 & 1.04e-05 & 1.16e-06  \\
256  & 64 $\times$ 64 $\times$ 1 & 1.05e-05 & 1.17e-06  \\
512  & 64 $\times$ 64 $\times$ 1 & 1.06e-05 & 1.20e-06  \\
1024 & 64 $\times$ 64 $\times$ 1 & 1.06e-05 & 1.23e-06   \\ \hline
\end{tabular}
\label{tab:wk_er}
\end{table}

\begin{table}[h!]
\centering
%\scriptsize
\caption{
A two-dimensional flux accuracy study for heterogeneous problem
with $4096$ cells. The problem represents a two-dimensional
slice of the heterogeneous weak scaling problem \eqref{table:ws.he.3d}
and \eqref{table:ws.he4.3d}
with $262$ thousand cells per subdomain.
}
\sizefonte
%\begin{tabularx}{\textwidth}{c|c|c|c|c|c}
\begin{tabular}{c|c|c|c}
%\hline \hline
%\multicolumn{5}{c}{\multirow{2}{*}{\textbf{Accuracy Study 2.2 (heterogeneous permeability)}}} \\
%\multicolumn{4}{c}{} \\
\hline \hline
\multicolumn{4}{c}{\multirow{2}{*}{\textbf{4096 cells per subd. (heterogeneous permeability)}}} \\
\multicolumn{4}{c}{} \\
\hline
\multicolumn{1}{c|}{} & Subd. cells & $\bar{H} = H$ & $\bar{H} = H/4$ \\
Cores  & ($n_x \times n_y \times n_z$) & $||\bu -\bu_h||_2/||\bu||_2$ & $||\bu -\bu_h||_2/||\bu||_2$  \\ \hline
32    & 64 $\times$ 64 $\times$ 1 & 7.91e-01 & 6.72e-01  \\
64    & 64 $\times$ 64 $\times$ 1 & 9.12e-01 & 8.77e-01 \\
128   & 64 $\times$ 64 $\times$ 1 & 8.54e-01 & 3.75e-01  \\
256   & 64 $\times$ 64 $\times$ 1 & 7.98e-01 & 2.67e-01  \\
512   & 64 $\times$ 64 $\times$ 1 & 6.20e-01 & 3.79e-01  \\
1024  & 64 $\times$ 64 $\times$ 1 & 6.24e-01 & 4.23e-01   \\ \hline
\end{tabular}                                                                         
\label{tab:wkhet_er}
\end{table}

%\begin{table}[h!]
%\centering
%%\scriptsize
%\caption{
%A two-dimensional flux accuracy study for heterogeneous problem
%with $4096$ cells. The problem represents a two-dimensional
%slice of the heterogeneous weak scaling problem \eqref{table:ws.he.3d}
%and \eqref{table:ws.he4.3d}
%with $262$ thousand cells per subdomain.
%}
%\sizefonte
%%\begin{tabularx}{\textwidth}{c|c|c|c|c|c}
%\begin{tabular}{c|c|c|c|c}
%%\hline \hline
%%\multicolumn{5}{c}{\multirow{2}{*}{\textbf{Accuracy Study 2.2 (heterogeneous permeability)}}} \\
%%\multicolumn{4}{c}{} \\
%\hline \hline
%\multicolumn{5}{c}{\multirow{2}{*}{\textbf{4096 cells per subd. (heterogeneous permeability)}}} \\
%\multicolumn{5}{c}{} \\
%\hline
%\multicolumn{1}{c|}{} & Subd. cells & $\bar{H} = H$ & $\bar{H} = H/4$ \\
%Cores  & ($n_x \times n_y \times n_z$) & $||\bu -\bu_h||_2/||\bu||_2$ & $||\bu -\bu_h||_2/||\bu||_2$ & Error Reduction \\ \hline
%32    & 64 $\times$ 64 $\times$ 1 & 7.91e-01 & 6.72e-01 & 11.9\%  \\
%64    & 64 $\times$ 64 $\times$ 1 & 9.12e-01 & 8.77e-01 & 3.5\%  \\
%128   & 64 $\times$ 64 $\times$ 1 & 8.54e-01 & 3.75e-01 & 47.9\%  \\
%256   & 64 $\times$ 64 $\times$ 1 & 7.98e-01 & 2.67e-01 & 53.1\%  \\
%512   & 64 $\times$ 64 $\times$ 1 & 6.20e-01 & 3.79e-01 & 24.1\%  \\
%1024  & 64 $\times$ 64 $\times$ 1 & 6.24e-01 & 4.23e-01 & 20.1\%   \\ \hline
%\end{tabular}                                                                         
%\label{tab:wkhet_er}
%\end{table}

%\clearpage

\section{Discussion}
To the best of our knowledge, there are very few 
weak and strong scaling results for parallel implementations of
multiscale methods in the
literature for three dimensional heterogeneous 
Darcy's flow and two-phase flow problems 
(see, \cite{manea16,manea17,manea19,ZYSW17,PED18}). 
In most of these papers the authors
consider a small number of processing cores (up to $20$), except in \cite{PED18}
that went up to $256$ cores in $16$ nodes, and numerical simulations with at
most hundreds of millions of cells which are comparable to our
smallest simulations. In addition, in these papers the most extensive
scalability study for multiscale methods is made for an algebraic extension of
multiscale methods \cite{manea16,manea17,manea19} that uses the multiscale
method as a preconditioner to damp low/high frequency modes of the resulting
discretized linear system related to the underlying elliptic PDE. 
It is noteworthy that in
this work we consider very large problems that can have billions of
cells, motivated by the numerical simulation of subsurface
flows, making use of a MPI base code for up to 1024 processing cores on 22 nodes. 
The numerical experiments reported in
the above mentioned papers are restricted to 128 million cells for both GPU
\cite{manea19} and CPU shared-memory architecture up to 20 processing cores on a single
node \cite{manea16,manea17}. 
In the study of \cite{PED18} that make use of $256$ processing cores, 
the three-dimensional simulations considered problems with about 16 million
discretized cells. 
Our results though, considered a larger set of nodes with an excellent
strong and weak scalability up to 1024 processing cores reaching up to $2$
billion cells in the weak scaling study and $1$ billion in the strong scaling
study.

%{\color{red}
%Moreover, it is important to note that although a strong scalability 
%up to 13 speedup for 20 processing cores is found in papers 
%\cite{manea16,manea17,manea19}, 
%our results considered a larger set of nodes with an excellent
%strong and weak scalability up to 1024 processing cores.
%}

\section{Concluding Remarks}\label{S:6} 

In this paper we developed a recursive formulation for the Multiscale Robin
Coupled Method that can be extended to the family of mixed
multiscale methods that the MRCM encompasses. The original 
global interface problem was
replaced by a set of small interface linear systems associated with adjacent
subdomains, in a hierarchy built as unions of nearest neighbor subdomains. 
A novel parallel algorithm is introduced and implemented for 
very large (up to 2 billions cells) problems, motivated
by the numerical simulation of subsurface flows. The recursive formulation was 
built to solve the global coarse algebraic problem more efficiently maintaining the
features of the underlying multiscale method. 
Through several numerical studies for both homogeneous and highly heterogeneous
permeability fields we showed that the new algorithm is very fast and exhibits
excellent scaling, with superlinear profile. As expected, the highly
heterogeneous problems present an increase of computational time compared to
the equivalent homogeneous problems. We observed small times for the 
numerical solution of the interface
problem, with the computation of the local boundary problems prevailing over 
the total computational time in all cases. 
Also, small changes on the intermediate coarse scale did not affect the
scalability of the formulation. The simulations were performed up to $1024$
processing cores without deterioration of the velocity field accuracy
presenting realistic potential of application in very large and highly heterogeneous
reservoirs. 

\section*{Acknowledgments}
The work presented here was partially funded by Petrobras research 
grants 2015/00398-0 and 2015/00400-4.
The authors also wish to thank
the Santos Dumont cluster located at the National Laboratory for 
Scientific Computing (LNCC) in Petr\'opolis, RJ, Brazil, and the Euler cluster 
at ICMC/USP in S\~ao Carlos, SP, Brazil.
E. Abreu was partially supported by CNPq 306385/2019-8 and PETROBRAS 2015/00398-0.

%E. Abreu thanks research grants as well as thanks to all the support given by
%the Brazilian funding agencies FAPESP 2019/20991-8 (São Paulo), CNPq 306385
%/2019-8 (National) and Petrobras 2015/00398-0 and 2019/00538-7.

%\section*{References}
\bibliography{paola_bib}

\end{document}